\documentclass[10pt,twocolumn]{IEEEtran}
\usepackage{amsmath,amssymb,amsthm,graphicx}
\usepackage{verbatim}
\graphicspath{{figures/}}
\usepackage{subfigure}
\usepackage[usenames]{color}
\definecolor{plum}  {rgb}{.4,0,.4}
\definecolor{BrickRed} {rgb}{0.6,0,0}
\usepackage[plainpages=false,
	pdfpagelabels,
	colorlinks=true,
	linkcolor=BrickRed,
	citecolor=plum]{hyperref}

\def\Reals{\mathbb{R}} 
\def\Naturals{\mathbb{N}} 

\def\sX{{\mathsf X}} 
\def\sU{{\mathsf U}} 

\def\cC{{\mathcal C}} 
\def\cF{{\mathcal F}} 
\def\cM{{\mathcal M}} 
\def\cP{{\mathcal P}} 

\def\PD{P^*} 
\def\DPop{\mathbb{T}} 

\def\E{\mathbb{E}} 
\def\twP{\check{P}} 
\def\twpi{\check{\pi}} 

\def\cT{{\mathcal T}} 

\def\wh#1{\widehat{#1}} 
\def\bd#1{\boldsymbol{#1}}
\def\argmin{\operatornamewithlimits{arg\,min}}
\def\supp{\operatorname{supp}}
\def\deq{\triangleq}

\newtheorem{theorem}{Theorem}
\newtheorem{lemma}{Lemma}
\newtheorem{proposition}{Proposition}

\newtheorem{assumption}{Assumption}
\newtheorem{remark}{Remark}

\begin{document}

\title{Online Markov Decision Processes\\
with Kullback--Leibler Control Cost}

\author{Peng Guan,~\IEEEmembership{Student Member,~IEEE}, Maxim Raginsky,~\IEEEmembership{Senior Member,~IEEE}, and~Rebecca~M.~Willett,~\IEEEmembership{Senior Member,~IEEE}
\thanks{This work was supported by NSF grant CCF-1017564 and by AFOSR grant FA9550-10-1-0390. A preliminary version of this work was presented at the American Control Conference, Montreal, Canada, June 2012.}
\thanks{P.~Guan is with the Department of Electrical and Computer Engineering, Duke University, Durham, NC 27708 USA (e-mail: peng.guan@duke.edu).}%
\thanks{M.~Raginsky is with the Department of Electrical and Computer Engineering and the Coordinated Science Laboratory, University of Illinois at Urbana-Champaign, Urbana, IL 61801 USA (e-mail: maxim@illinois.edu).}%
\thanks{R.~Willett is with the Department of Electrical and Computer Engineering, University of Wisconsin-Madison, Madison, WI 53796 USA (e-mail: willett@wisc.edu).
}%
}

\markboth{}{}

\maketitle

\begin{abstract}
This paper considers an online (real-time) control problem that
involves an agent performing a discrete-time random walk over a finite
state space. The agent's action at each time step is to specify the
probability distribution for the next state given the current state.
Following the set-up of Todorov, the state-action cost at each time
step is a sum of a state cost and a control cost given by the
Kullback-Leibler (KL) divergence between the agent's next-state
distribution and that determined by some fixed passive dynamics. The
online aspect of the problem is due to the fact that the state cost
functions are generated by a dynamic environment, and the agent learns
the current state cost only after selecting an action. An explicit
construction of a computationally efficient strategy with small regret
(i.e., expected difference between its actual total cost and the
smallest cost attainable using noncausal knowledge of the state costs)
under mild regularity conditions is presented, along with a
demonstration of the performance of the proposed strategy on a
simulated target tracking problem. A number of new results on Markov
decision processes with KL control cost are also obtained.
\end{abstract}


\IEEEpeerreviewmaketitle

\thispagestyle{empty}

\section{Introduction}

Markov decision processes (MDPs)
\cite{Puterman,LermaLasserreMDP,AC_MDP_survey} comprise a popular framework
for sequential decision-making in a random dynamic environment. At
each time step, an agent observes the state of the system of interest
and chooses an action. The system then transitions to its next state,
with the transition probability determined by the current state and
the action taken. There is a (possibly time-varying) cost associated
with each admissible state-action pair, and a policy (feedback law)
for mapping states to actions is selected to minimize average cost. In
the basic MDP framework, it is assumed that the cost functions and the
transition probabilities are known in advance, the policy is designed
``offline'' (e.g., using dynamic programming), and the optimality
criterion is forward-looking, taking into account the effect of past
actions on future costs. In many practical problems, however, this
degree of advance knowledge is unavailable. When neither the transition probability nor the cost functions are known in advance, various reinforcement learning (RL) methods, such as the celebrated $Q$-learning algorithm \cite{Watkins_Dayan_QLearning,Tsitsiklis_QLearning} and its variants, can be used to learn an optimal policy in an online regime. However, the key assumptions underlying RL are that the agent is operating in a stochastically stable environment, and that the state-action costs (or at least their expected values with respect to any environmental randomness) do not vary with time. These assumptions are needed to ensure that the agent is eventually able to learn an optimal stationary control policy.

Another framework for sequential decision-making, dating back to the seminal work of Robbins \cite{Robbins_compound} and Hannan \cite{Hannan} and now widely used in the
machine learning community \cite{PLG}, deals with nonstochastic, unpredictable environments. In this {\em online learning} (or {\em sequential prediction}) framework, the effects of the
environment are modeled by an arbitrarily varying sequence of cost functions,
where the cost function at each time step is revealed to the agent
only {\em after} an action has been taken. There is no
state, and the goal of the agent is to minimize {\em regret}, i.e.,
the difference between the total cost incurred using causally
available information and the total cost of the best single action
that could have been chosen in hindsight. In contrast with MDPs, the
regret-based optimality criterion is necessarily myopic and
backward-looking, since the cost incurred at each time step depends
only on the action taken at that time step, so past actions have no
effect on future costs. There is also a more stringent model of online learning, in which the agent observes not the entire cost function for each time step, but only the value of this cost at the currently taken action \cite{Auer_bandits}. This model is inspired by the celebrated {\em multiarmed bandit} problem first introduced by Robbins \cite{Robbins}, and is referred to as the {\em nonstochastic bandit problem}. One widely used way of constructing regret-minimizing strategies for such bandit problems is to randomize the agent's actions ({\em exploration}) so that the random cost value revealed to the agent can be used to construct an unbiased estimate of the full cost function, which is then fed into a suitable strategy that minimizes regret under the assumption of full information ({\em exploitation}). We will not consider nonstochastic bandit problems in this paper. Instead, we refer the reader to a recent survey by Bubeck and Cesa-Bianchi \cite{Bubeck_bandits} that discusses both stochastic and nonstochastic bandit problems.

Recent work by Even-Dar et al. \cite{EvenDar} and Yu et al. \cite{Yu} combines the MDP and the online learning frameworks into what may be described as {\em online MDPs} with finite state and action spaces. Like in the traditional MDP setting, the agent observes the current state and chooses an action, and the system transitions to the next state according to a fixed and known Markov law. However, like in the online framework, the one-step cost functions form an arbitrarily varying sequence, and the cost function corresponding to each time step is revealed to the agent after the action has been taken. The objective of the agent is to minimize regret relative to the best stationary Markov policy that could have been selected with full knowledge of the cost function sequence over the horizon of interest. The time-varying cost functions may represent unmodeled aspects of the environment or collective (and possibly irrational) behavior of any other agents that may be present; the regret minimization viewpoint then ensures that the agent's {\em online} policy is robust against these effects.

\subsection{Brief problem statement and motivating examples}
\label{ssec:setup_examples}

We give here a brief statement of the problem of interest in order to fix ideas; a more detailed formulation is given later on. The reader may wish to consult Section~\ref{ssec:notation} for notation.

The set-up considered in \cite{EvenDar,Yu} is motivated by problems in machine learning and artificial intelligence, where the actions are the main object of interest, and the state merely represents memory effects present in the system. In this paper, we take a more control-oriented view: the emphasis is on steering the system along a desirable state trajectory through actions selected according to a state feedback law. Following the formulation proposed recently by Todorov \cite{TodorovNIPS,TodorovCDC,TodorovPNAS}, we allow the agent to modulate the state transitions directly, so that actions (resp., state feedback laws) correspond to probability distributions (resp., Markov kernels) on the underlying state space. As in \cite{TodorovNIPS,TodorovCDC,TodorovPNAS}, the one-step cost is a sum of two terms: the state cost, which measures how ``desirable'' each state is, and the control cost, which measures the deviation of the transition probabilities specified by the chosen action from some fixed {\em default} or {\em passive dynamics}. (We also refer the reader to a recent paper by Kappen et al.~\cite{Kappen_etal}, which interprets Todorov's set-up as an inference problem for probabilistic graphical models.)

More precisely, we consider an MDP with a finite state space $\sX$, where the action space $\sU$ is the simplex $\cP(\sX)$ of probability distributions over $\sX$. A fixed Markov matrix (transition kernel) $\PD = [\PD(x,y)]_{x,y \in \sX}$ is given. A stationary Markov policy (state feedback law) is a mapping $w : \sX \to \cP(\sX)$, so if the system is in state $x \in \sX$, then the transition to the next state is stochastic, as determined by the probability distribution $u(\cdot) = w(x) \in \cP(\sX)$. In other words, if we denote the next state by $X^+$, then the state transitions induced by the action $u$
are governed by the conditional probability law
$$
\Pr \{X^+ = x^+ | X=x \} = P(x, x^+) = u(x^+) = [w(x)](x^+).
$$
The one-step state-action cost $c(x,u)$ consists of two terms, the
{\em state cost} $f(x)$, where $f : \sX \mapsto \Reals_+$ 
is a given function,
and the {\em control cost}, which penalizes any deviation of the
next-state distribution $u(\cdot) = w(x)$ from the one prescribed by
$\PD(x,\cdot)$, the row of $\PD$ corresponding to $x$. To
motivate the introduction of such control costs, we can imagine the
situation, in which implementing the state transitions according to
$\PD$ can be done ``for free.'' However, it may very well be the case
that following $\PD$ will be in conflict with the goal of keeping the
state cost low. From this perspective, it may actually be desirable to
deviate from $\PD$. Any such deviation may be viewed as an {\em active
  perturbation} of the {\em passive dynamics} prescribed by $\PD$, and
the agent should attempt to balance the tendency to keep the state
costs low against allowing too strong of a perturbation of $\PD$. Our
choice of control cost is inspired by the work of Todorov
\cite{TodorovNIPS,TodorovPNAS}, and is given by the {\em
  Kullback--Leibler divergence} (or the {\em relative entropy})
\cite{CoverThomas} $D(u\|\PD(x,\cdot))$ between the proposed next-state
distribution $u(\cdot)$ and the next-state distribution prescribed by
the passive dynamics $\PD$. One useful property of this control cost
is that it automatically forbids all those state transitions that are
already forbidden by $\PD$. Indeed, if for a given $x \in \sX$ there exists some $y \in \sX$
such that $u(y) = [w(x)](y) > 0$, while $\PD(x,y) = 0$, then
$D(u\| \PD(x,\cdot)) = + \infty$. Thus, the overall
one-step state-action cost is given by
\begin{align}\label{eq:Todorov_cost}
c(x,u) = f(x) + D(u \| \PD(x,\cdot)), \quad \forall x \in \sX, u \in \cP(\sX).
\end{align}
In the online version of this problem (detailed in Section~\ref{ssec:model}), the state costs form an arbitrarily varying sequence $\{f_t\}^\infty_{t=1}$, and the agent learns the state cost for each time step only after having selected the transition law to determine the next state. For any given value of the horizon, the regret is computed with respect to the best stationary Markov policy (state feedback law) that could have been chosen in hindsight. The precise definition of regret is given in Section~\ref{ssec:regret}. 

Since this is a nonstandard set-up, we take a moment to situate it in the context of usual models of MDPs. In a standard MDP with finite state and action spaces, we have a finite collection of Markov matrices $P_u$ on $\sX$ indexed by the actions $u \in \sU$. State feedback laws are functions $w : \sX \to \sU$, and the set of all such functions is finite with cardinality $ |\sU|^{|\sX|}$. Therefore, in each state $x \in \sX$ the agent has at most $|\sU|^{|\sX|}$  choices for the distribution of the next state $X^+$, and we may equivalently represent each state feedback law $w$ as a mapping from $\sX$ into $\cP(\sX)$ with $x \mapsto P_{w(x)}(x,\cdot)$. Since the state space $\sU$ is finite, the range of this mapping is a finite subset of the probability simplex $\cP(\sX)$. The criterion for selecting this next-state distribution pertains to minimization of the expectation of the immediate state-action cost plus a suitable value function that accounts for the effect of the current action on future costs. In many cases, the one-step state-action cost $c(x,u)$ decomposes into a sum of state cost $f(x)$ and control cost $g(x,u)$, where $f(x)$ quantifies the (un)desirability of the state $x$, while $g(x,u)$ represents the effort required to apply action $u$ in state $x$.

In the set-up of \cite{TodorovNIPS,TodorovCDC,TodorovPNAS},
  the collection of all possible next-state distributions is
  unrestricted. As a consequence, any mapping $w : \sX \to \cP(\sX)$
  is a feasible stationary Markov policy. Since any Markov matrix $P$
  on $\sX$ can be equivalently represented as a mapping from $\sX$
  into $\cP(\sX)$ with $x \mapsto P(x,\cdot)$, there is thus a
  one-to-one correspondence between state feedback laws and Markov
  matrices on $\sX$. In contrast to the case when the agent may choose
  among a finite set of actions, the probability simplex $\cP(\sX)$ is
  an uncountable set, so the agent has considerably greater freedom to choose the next state distribution. As before, we introduce a state-action cost $c(x,u) = f(x) + g(x,u)$, where $f(x)$ measures the (un)desirability of state $x$, while $g(x,u)$ quantifies the difficulty of executing action $u$ in state $x$. Since actions $u$ now correspond to probability distributions, and we choose the Kullback--Leibler control cost $g(x,u) = D(u \| \PD(x,\cdot))$, where $\PD$ is a fixed Markov matrix on the state space $\sX$ that may represent, e.g., the ``free'' dynamics of the system in the absence of external controls.

The Kullback--Leibler divergence is widely used in stochastic control and inference. First of all, it has many desirable properties, such as nonnegativity and convexity \cite{CoverThomas}. Secondly, if we adopt the viewpoint that the purpose of a control policy is to shape the joint distribution of all relevant variables describing the closed-loop behavior of the system, then using the relative entropy to compare the distribution induced by any control law to some reference model leads to functional equations for the optimal policy that are often easier to solve than the corresponding dynamical programming recursion \cite{Karny,Sindelar_etal} (e.g., see \cite{Karny} for an alternative derivation of the optimal controller in an LQG problem using relative entropy instead of dynamic programming); similar ideas are fruitful in the context of robust control, where the relative entropy is used to quantify the radius of uncertainty around some nominal system \cite{Petersen_etal,Charalambous,Robustness_book}. Moreover, the relative entropy is a canonical regularization functional for stochastic nonlinear filtering problems \cite{Mitter}: an optimal Bayesian filter is the solution of a variational problem that entails minimization of the sum of expected negative log-likelihood (which can be interpreted as state cost) and a relative entropy with respect to the prior measure on the state space.

To further motivate our interest in problems of this sort, let us consider two
examples. One is {\em target tracking} with an arbitrarily moving
target (or multiple targets). In this example, the state space $\sX$
is the vertex set of an undirected graph, and the passive dynamics
$\PD$ specifies some default {\em random walk} on this graph. The
tracker's discrete-time motion is constrained by the topology of
the graph, while the targets' motions are not. At each time $t$, the state cost $f_t$ is the tracking error, which quantifies how far the tracker is from the targets. For instance, it may be given by the graph distance (length of shortest path) between the tracker's current location and the location of the closest target. Other possibilities can also be considered, including some based on noisy information on the location of the targets. The control cost penalizes the tracker's deviation from $\PD$ as it attempts to track the targets. The passive dynamics $\PD$ can be seen as the tracker's prior model for the targets' motion. Moreover, if $\PD$ is sufficiently rapidly mixing, then any tracker that follows $\PD$ will visit every vertex of the graph infinitely often with probability one; however, there is no guarantee that the tracker's prior model is correct (i.e., that the tracker will be anywhere near the targets). Hence, the state-action cost will trade off the tendency of the tracker to ``cover'' the graph as much as possible (exploration) against the tendency to follow a potentially faulty model of the targets (exploitation). 

Another example setting is real-time control of a {\em brain-machine
  interface}. There, the state space $\sX$ may be the set of possible
positions or modes of a neural prosthetic device, and the passive
dynamics $\PD$ may encode the ``natural'' (free) dynamics of the
device in the absence of user control; we may assume, for instance,
that the state transitions prescribed by $\PD$ correspond to
``minimum-energy'' operating mode of the device. If the user wishes to
make the device execute some trajectory, the state cost $f_t$ at time
$t$ may represent the deviation of the current point on the trajectory
from the one intended by the user. Since the user is a human operator
with conscious intent, we may not want to ascribe an a priori model to
her intended trajectory, and instead treat it as an individual sequence modulating the state costs $\{f_t\}^\infty_{t=1}$. In this setting, the Kullback--Leibler control cost penalizes significant deviations from the free dynamics $\PD$, since these will typically be associated with energy expenditures.

The common thread running through these two examples (and it is certainly possible to construct many others) is that they model real-time interaction of a particular system with some well-defined ``reference'' or ``nominal" dynamics $\PD$ with a potentially unpredictable environment (which may include hard-to-model adversaries or rational agents, etc.), and we must balance the tendency to respond to immediate changes in the environment against the need to operate the system near the nominal mode. Since no offline policy design is possible in such circumstances, the regret minimization framework offers a meaningful alternative.

\subsection{Our contributions and comparison with relevant literature}
\label{ssec:contribution}

In this paper, we give an explicit construction of a strategy for the agent, such that the regret relative to any uniformly ergodic class of stationary Markov policies grows {\em sublinearly} as a function of the horizon. The only regularity conditions needed for this result to hold are (a) uniform boundedness of the state costs (the agent need not know the bound, only that it exists); and (b) ergodicity of the passive dynamics. Moreover, our strategy is computationally efficient: the time is divided into phases of increasing length, and during each phase the agent applies a stationary Markov policy optimized for the average of the state cost functions revealed during all of the preceding phases. Thus, our strategy belongs to the class of so-called ``lazy'' strategies for online decision-making problems \cite{Vovk,Merhav_etal,Kalai_Vempala}; a similar approach was also taken by Yu et al.~\cite{Yu} in their paper on online MDPs with finite state and action spaces. The main advantage of lazy strategies is their computational efficiency, which, however, comes at the price of suboptimal scaling of the regret with the time horizon. We comment on this issue further in the sequel.

Our main contribution is an extension of the theory of online MDPs to a wide class of control problems that lie outside the scope of existing approaches \cite{EvenDar,Yu}. More specifically:
\begin{enumerate}
	\item While in \cite{EvenDar,Yu} both the state and the action spaces are finite, we only assume this for the state space. Our action space is the simplex of probability distributions on the state space, which is a compact subset of a Euclidean space. Hence, the techniques used in the existing literature are no longer directly applicable. (It is also possible to extend our approach to continuous state spaces, but additional regularity conditions will be needed. This extension will be the focus of our future work.)
	\item Yu et al. \cite{Yu} assume that the underlying MDP is unichain \cite[Sec.~8.3]{Puterman} and satisfies a certain uniform ergodicity condition (a similar assumption is also needed by Even-Dar et al. \cite{EvenDar}). Their assumption is rather strong, since it places significant simultaneous restrictions on an exponentially large family of Markov chains on the state space (each chain corresponds to a particular choice of state feedback law, and there are $|\sU|^{|\sX|}$ such laws). It is also difficult to verify, since the problem of determining whether an MDP is unichain is NP-hard \cite{Tsitsiklis}. By contrast, our ergodicity assumption pertains to only {\em one} Markov chain (the passive dynamics $\PD$), it can be efficiently verified in polynomial time, and we prove that it automatically implies uniform ergodicity of all stationary control laws that could possibly be invoked by our strategy.
	\item Because these stationary control laws correspond to solutions of certain average-cost optimality equations (ACOEs) in the set-up of Todorov \cite{TodorovNIPS,TodorovCDC,TodorovPNAS}, we establish and subsequently exploit several useful and previously unknown results concerning the continuity and uniform ergodicity of optimal policies for Todorov's problem. These results, as well as the techniques used to prove them, play a very important role in our overall contribution. Indeed, in the online setting, the state cost functions are revealed to the agent in real time. Hence, any policy used by the agent must rely on estimates (or forecasts) of future state costs based on currently available information. Our new results on Todorov's optimal control laws provide sharp bounds on the sensitivity of these laws to misspecification of state costs, and may be of independent interest.
	\item In \cite{Yu}, the policy computation at the beginning of each phase requires solving a linear program and then adding a carefully tuned random perturbation to the solution. As a result, the performance analysis in \cite{Yu} is rather lengthy and technical (in particular, it invokes several advanced results from perturbation theory for linear programs). By contrast, even though we are working with a continuous action space, all policy computations in our case reduce to solving finite-dimensional eigenvalue problems, without any need for additional randomization. Moreover, even though the overall scheme of our analysis is similar to the one in \cite{Yu} (which, in turn, is inspired by existing work on lazy strategies \cite{Vovk,Kalai_Vempala,Merhav_etal}), the proof is self-contained and much less technical, relying on our new results pertaining to Todorov-type optimal control laws.		
\end{enumerate}
A preliminary version of this work has appeared in a conference publication \cite{MaxBeccaPeng}, and most of the proofs were omitted due to space limitations. Since a major part of our contribution is a set of probabilistic analysis techniques for MDPs with Kullback--Leibler control cost (in both online and offline settings), the present paper fills in the missing details. In addition, most of our new results on the sensitivity of Todorov-type optimal controllers to perturbations of state costs were omitted from \cite{MaxBeccaPeng}. The present paper not only gives a self-contained treatment of these results, but also demonstrates their crucial role in performance analysis of online strategies for Todorov-type MDPs. Finally, compared to \cite{MaxBeccaPeng}, the present paper reports a more thorough and improved empirical evaluation of our proposed strategy in the context of target tracking on a large graph. In particular, we report the results of Monte-Carlo simulation of our strategy (with error bars) and compare it to two baseline strategies: (a) the best stationary policy that could be chosen with full prior knowledge of the state cost sequence and (b) the best stationary policy chosen from a large pool of randomly sampled policies {\em without} advance knowledge of state costs. In \cite{MaxBeccaPeng}, we only compared our strategy to the passive dynamics $\PD$. The experimental results reported here show that (a) the regret of our strategy w.r.t.\ the best stationary policy chosen in hindsight is nonnegative and grows sublinearly with time (thus validating our theoretical bound), and (b) in simulations, our strategy performs strictly better than any randomly sampled stationary policy.

\subsection{Organization of the paper}

The remainder of the paper is organized as follows. We close this section with a brief summary of frequently used notation. Section~\ref{sec:problem} contains precise formulation of the online MDP problem and presents our main result, Theorem~\ref{thm:main_r}. In preparation for the proof of the theorem, Section~\ref{sec:prelims} contains preliminaries on MDPs with KL control cost \cite{TodorovNIPS,TodorovCDC,TodorovPNAS}, including a number of new results pertaining to optimal policies.  Section~\ref{sec:strategy} then describes our proposed strategy, whose performance is then analyzed in Section~\ref{sec:proof} in order to prove Theorem~\ref{thm:main_r}. Some simulation results are presented in Section~\ref{sec:simulations}. We close by summarizing our contributions and outlining some directions for future work. Proofs of all intermediate results are relegated to the Appendix.

\subsection{Notation}
\label{ssec:notation}

We will denote the underlying finite state space by $\sX$. A matrix $P = [P(x,y)]_{x,y \in \sX}$ with nonnegative entries, and with the rows and the columns indexed by the elements of $\sX$, is called {\em stochastic} (or {\em Markov}) if its rows sum to one: $\sum_{y \in \sX} P(x,y) = 1, \forall x \in \sX$.

We will denote the set of all such stochastic matrices  by $\cM(\sX)$, the set of all probability distributions over $\sX$ by $\cP(\sX)$, the set of all functions $f : \sX \to \Reals$ by $\cC(\sX)$, and the cone of all nonnegative functions $f : \sX \to \Reals_+$ by $\cC_+(\sX)$. We will represent the elements of $\cP(\sX)$ by row vectors and denote them by $\pi,\mu,\nu$, etc., and the elements of $\cC(\sX)$ by column vectors and denote them by $f,g,h$, etc. The total variation (or $L_1$) distance between $\mu, \nu \in \cP(\sX)$ is
\begin{align*}
\| \mu - \nu \|_1 \deq \sum_{x \in \sX} \lvert \mu(x) - \nu(x) \rvert.
\end{align*}
The {\em Kullback--Leibler divergence} (or {\em relative entropy}) \cite{CoverThomas} between $\mu$ and $\nu$ is
\begin{displaymath}
D(\mu \| \nu) \deq  \begin{cases}
\sum_{x \in \sX} \mu(x) \log \dfrac{\mu(x)}{\nu(x)} & \textrm{if ${\rm supp}(\mu) \subseteq {\rm supp}(\nu$)} \\
+ \infty & \textrm{otherwise}
\end{cases}
\end{displaymath}
where ${\rm supp}(\mu) \deq \{ x \in \sX: \mu(x) > 0 \}$ is the {\em support} of $\mu$. Here and in the sequel, we work with natural logarithms. The {\em span seminorm} (also called the {\em oscillation}) of $f \in \cC(\sX)$ is defined as
\begin{align*}
\| f \|_s \deq \max_{x \in \sX} f(x) - \min_{x \in \sX} f(x).
\end{align*}
Note that $\| f \|_s = 0$ if and only if $f(x) = c$ for some constant $c \in \Reals$ and all $x \in \sX$; $\| f \|_s = \| f + c \|_s$ for any $f \in \cC(\sX)$ and $c \in \Reals$. We also define the {\em sup norm} $\| f \|_\infty \deq \max_{x \in \sX} \lvert f(x) \rvert$ and note that $\| f \|_s \le 2 \| f \|_\infty$.

Any Markov matrix $P \in \cM(\sX)$ acts on probability distributions from the right and on functions from the left:
\begin{align*}
	\mu P(y) = \sum_{x \in \sX}\mu(x)P(x,y), \qquad Pf(x) = \sum_{y \in \sX}P(x,y)f(y).
\end{align*}
We say that $P$ is {\em unichain} \cite{LermaLasserreMC} if the corresponding Markov chain has a single recurrent class of states (plus a possibly empty transient class). The is equivalent to $P$ having a unique invariant distribution $\pi_P$ (i.e. $\pi_P P = \pi_P$) \cite{Seneta}. We will denote the set of all such Markov matrices over $\sX$ by $\cM_1(\sX)$. Given $\rho \in [0, 1]$, we say that $P$ is $\rho$-{\em contractive} if
\begin{align*}
\| \mu P - \nu P \|_1 \le \rho \| \mu - \nu \|_1, \qquad \forall \mu, \nu \in \cP(\sX)
\end{align*}
(in fact, every $P \in \cM(\sX)$ is 1-contractive). We will denote the set of $\rho$-contractive Markov matrices by $\cM_1^{\rho}(\sX)$. It is easy to show that, for every $0 \le \rho <1$, $\cM_1^{\rho}(\sX) \subset \cM_1(\sX)$. The  {\em Dobrushin ergodicity coefficient} \cite{Seneta,Cappe} of $P \in  \cM(\sX)$ is given by 
\begin{align*}
	\alpha(P) \deq \frac{1}{2} \max_{x,x' \in \sX} \| P(x,\cdot) - P(x',\cdot) \|_1,
\end{align*}
and it can be shown that any $P \in \cM(\sX)$ is $\alpha(P)$-contractive \cite{Seneta,Cappe}. Finally, for any $P, P' \in  \cM(\sX)$ we define the supremum distance
\begin{align*}
	\| P - P' \|_\infty \deq \max_{x \in \sX} \| P(x,\cdot) - P'(x,\cdot) \|_1.
\end{align*}

\section{Problem formulation and the main result}
\label{sec:problem}

\subsection{The model}
\label{ssec:model}

Given the finite state space $\sX$, let $\cF$ be a fixed subset of $\cC_+(\sX)$, and let $x_1 \in \sX$ be a fixed initial state. Consider an agent (A) performing a controlled random walk on $\sX$ in response to a dynamic environment (E). The interaction between A and E proceeds as follows:
\begin{center}
\begin{tabular}{|l|}
\hline
$X_1=x_1$\\
for $t = 1,2,\ldots$\\
\ \ A selects $P_t \in \cM(\sX)$ and draws $X_{t+1} \sim P_t(X_t,\cdot)$\\
\ \ E selects $f_t \in \cF$ and announces it to A\\
end for\\
\hline 
\end{tabular}
\end{center}
At each $t \ge 1$, A selects the transition probabilities $P_t(x,y) = \Pr\{X_{t+1}=y|X_t=x\}$ based on his knowledge $f^{t-1} = (f_1,\ldots,f_{t-1})$, and incurs the {\em state cost} $f_{t}(X_{t})$ and the {\em control cost} $D(P_t(X_t,\cdot) \| \PD(X_t,\cdot))$. The total cost incurred by the agent A at time $t$ is given by
\begin{align*}
	c_t(X_t,P_t) = f_t(X_t) + D(P_t(X_t,\cdot)\| \PD(X_t,\cdot)).
\end{align*}
and the objective is to minimize a suitable notion of regret. 

\subsection{Strategies and regret}
\label{ssec:regret}

A {\em strategy} for the agent A is a collection of mappings $\gamma = \lbrace \gamma_t \rbrace_{t=1}^{\infty}$ where $\gamma_t : \mathcal F^{t-1} \rightarrow \cM(\sX)$, so that $P_t = \gamma_t(f^{t-1})$. This means our strategy is based on the complete knowledge of all the past cost functions. The cumulative cost of $\gamma$ after $T$ steps is
\begin{align*}
C_T = \displaystyle \sum_{t=1}^T c_t(X_t, P_t) = \displaystyle \sum_{t=1}^T c_t(X_t, \gamma_t(f^{t-1})).
\end{align*}
To define the regret after $T$ steps, we will consider the gap between $C_T$ and the expected cumulative cost that A could have achieved in hindsight by using a stationary unichain random walk on $\sX$ (with full knowledge of $f^T$). This gap arises through the agent's lack of prior knowledge on the sequence of state cost functions. Formally, we define the regret of $\gamma$ after $T$ steps w.r.t.\ a particular $P \in \cM_1(\sX)$ by\footnote{To keep the notation clean, we will suppress the dependence of the cumulative cost $C_T$ and the regret $R_T$ on the strategy $\gamma$ and on the state costs $f_1,\ldots,f_T$.}
\begin{align*}
	R_T(P) \deq C_T - \E_{x_1}^P \left[  \sum_{t=1}^{T} c_t(X_t, P) \right],
\end{align*}
where the expectation is taken over the Markov chain induced the by {\em comparison} transition kernel $P$ with initial state $X_1 = x_1$. In this work, we make the following basic assumption concerning the environment E:
\setcounter{assumption}{-1}
\begin{assumption}[Oblivious environment]\label{as:oblivious} The environment {\em E} is {\em oblivious} (or {\em nonadversarial}), i.e., for every $t$, $f_t$ depends only on $f^{t-1}$, but not on $X^t$. \end{assumption}
\noindent Assumption~\ref{as:oblivious} is standard in the literature on sequential prediction \cite{PLG} (in particular, it is also imposed by Yu et al.\ \cite{Yu}). In our case, it implies that, for a fixed sequence $f_1,f_2,\ldots$ of state costs chosen by E, the state process $\bd{X} = \{X_t\}^\infty_{t=1}$ induced by A's choices $P_1,P_2,\ldots$ is a (time-inhomogeneous) Markov chain. Now consider some set ${\mathcal N} \subset \cM_1(\sX)$. Adopting standard terminology \cite{PLG}, we will say that $\gamma$ is {\em Hannan-consistent} w.r.t.\ $\mathcal N$ if
 \begin{align}\label{eq:hannan}
	\limsup_{T \to \infty} \sup_{P \in {\mathcal N}} \sup_{f_1, . . . , f_T \in \cF}  \frac{\E R_T(P)}{T} \le 0, 
\end{align}
where the expectation is w.r.t.\ the law of the process $\bd{X}$ starting at $X_1 = x_1$. In other words, a strategy is Hannan-consistent if its worst case (over $\cF$) expected per-round regret converges to zero uniformly over $\mathcal N$. While it is certainly true that some nonstationary policy with complete prior knowledge of the state cost sequence may (and will) outperform any stationary policy, we limit our consideration to stationary reference policies for two reasons. One is the need to have a fair comparison: indeed, no truly online strategy could compete with the best (i.e., omniscient) nonstationary policy. The other is that we can alternatively interpret the Hannan consistency condition \eqref{eq:hannan} as follows: as the horizon $T$ increases, the smallest average cost achievable by a strategy which is Hannan-consistent w.r.t.\ ${\mathcal N}$ will converge to the smallest long-term average cost achievable by any stationary Markov strategy in ${\mathcal N}$ on an MDP with the state cost given by the {\em empirical average} $(1/T)\sum^T_{t=1} f_t$ of the state costs revealed up to time $T$.

\subsection{The main result}
\label{ssec:main_result}

Our main result (Theorem~\ref{thm:main_r} below) guarantees the existence of a Hannan-consistent strategy against any uniformly ergodic collection of stationary unichain policies under the following two assumptions on the passive dynamics $\PD$:

\begin{assumption}[Irreducibility and aperiodicity]\label{as:irred} The passive dynamics $\PD$ is {\em irreducible} and {\em aperiodic}, where the former means that, for every $x,y \in \sX$, there exists some $n \in \Naturals$ such that $(\PD)^n(x,y) > 0$, while the latter means that, for every $x \in \sX$, the greatest common divisor of the set $\{ n \in \Naturals: (\PD)^n(x,x) > 0 \}$ is equal to $1$.
\end{assumption}
\begin{assumption}[Ergodicity]\label{as:Dobrushin} The Dobrushin ergodicity coefficient $\alpha(\PD)$ is strictly less than $1$.
\end{assumption}
\noindent Assumption~\ref{as:irred} ensures that $\PD$ has a unique everywhere positive invariant distribution $\pi^*$ \cite{Seneta} and, for a finite $\sX$, it is equivalent to the existence of some $\bar{n} \in \Naturals$, such that
\begin{align*}
	\theta \deq \min_{x,y \in \sX} (\PD)^{\bar{n}}(x,y) > 0
\end{align*}
(see, e.g., Theorem~1.4 in \cite{Seneta}). Assumption~\ref{as:Dobrushin}, which is frequently used in the study of MDPs with average cost criterion \cite{LermaAdaptive,MasiStettner,AC_MDP_survey}, guarantees that the convergence to $\pi^*$ is exponentially fast (so that $\PD$ is geometrically ergodic), and it also imposes a stronger type of ergodicity, since a Markov matrix $P \in \cM(\sX)$ has $\alpha(P) < 1$ if and only if for any pair $x, x' \in \sX$ there exists at least one $y \in \sX$, such that $y$ can be reached from both $x$ and $x'$ in one step with strictly positive probability. For example, if $\PD$ satisfies the {\em Doeblin minorization condition} \cite{Cappe,MeynTweedie}, i.e., if there exist some $\delta \in (0,1]$ and some $\mu \in \cP(\sX)$, such that $\PD(x,y) \ge \delta \mu(y)$ for all $x,y \in \sX$, then we will have $\alpha(\PD) \le 1-\delta < 1$ (see, e.g., Lemma~4.3.13 in \cite{Cappe}). 

\begin{remark}{\em As we pointed out in Section~\ref{ssec:contribution}, our assumption is actually much milder than the assumptions made in related literature \cite{EvenDar,Yu,NeuCsaba}. Recent work by Neu et al. \cite{NeuCsaba} shows that the assumptions made in \cite{EvenDar,Yu,NeuCsaba} are valid only if the Dobrushin ergodicity coefficient of the state transition kernel induced by every policy is strictly smaller than one. By contrast, we only assume this for the passive dynamics $\PD$.}
\end{remark} 

With these assumptions in place, we are now ready to state our main result:
\begin{theorem}\label{thm:main_r} Let $\mathcal F$ consist of all $f \in \cC_+(\sX)$ with $\| f \|_\infty \le 1$. Fix an arbitrary $\epsilon \in (0,1/3)$. Under Assumptions~\ref{as:oblivious}--\ref{as:Dobrushin}, there exists a strategy $\gamma$, such that for any $\rho \in [0,1)$, 
\begin{align}\label{eq:thm1}
	\sup_{P \in \cM_1^{\rho}(\sX)}\sup_{f_1, . . . , f_T \in \mathcal F} \frac{\E R_T(P)}{T} = O(T^{-1/4+\epsilon}).
\end{align}
As a consequence, the strategy $\gamma$ is Hannan-consistent w.r.t.\ $\cM_1^{\rho}(\sX)$.
\end{theorem}

\begin{remark}{\em The constant hidden in the $O(\cdot)$ notation depends only on the passive dynamics $\PD$ and on the contraction rate $\rho$ of the baseline policies in $\cM_1^{\rho}(\sX)$; cf.~Eq.~\eqref{eq:thm2}, and the discussion preceding it, for details.}
\end{remark} 

\section{Preliminaries}
\label{sec:prelims}

Our construction of a Hannan-consistent strategy in Theorem~\ref{thm:main_r} relies on Todorov's theory of MDPs with KL control cost \cite{TodorovNIPS, TodorovCDC,TodorovPNAS}. In this section, we give an overview of this theory and present several new results that will be used later on.

First, let us recall the general set-up for MDPs with finite state space $\sX$ and compact action space $\sU$ under the average cost criterion (see e.g., \cite{LermaLasserreMDP} or \cite{AC_MDP_survey}). It involves a family of Markov matrices $P_u \in \cM(\sX)$ indexed by actions $u \in \sU$. The (long-term) \emph{average cost} of a stationary Markov policy (state feedback law) $w: \sX \rightarrow \sU$ with initial state $X_1 = x_1$ is given by
\begin{align}\label{eq:long_avr}
J(w,x_1) \deq \limsup_{T \rightarrow \infty} \frac{1}{T} \E_{x_1}^w \left[ \sum_{t=1}^T c(X_t, w(X_t)) \right],
\end{align} 
where the expectation $\E_{x_1}^w[\cdot]$ is w.r.t.\ the law of the Markov chain $\bd{\mathit X} = \lbrace X_t \rbrace$ with controlled transition probabilities
\begin{align*}
\Pr\{ X_{t+1} = y \vert X_t = x \} = P_{w(x)}(x,y), \qquad X_1 = x_1
\end{align*}
and $c: \sX \times \sU \rightarrow \Reals$ is the one step state-action cost. The construction of an optimal policy to minimize \eqref{eq:long_avr} for every  $x_1$ revolves around the {\em average-cost optimality equation} (ACOE)
\begin{align}\label{eq:ACOE}
	h(x) + \lambda = \min_{u \in \sU(x)} \left\{ c(x,u) + P_u h(x)\right\}, \,\, x \in \sX
\end{align}
where $\sU(x) \subseteq \sU$ is the set of allowable actions in state $x$. If a solution pair $(\lambda,h) \in \Reals \times \cC(\sX)$ exists with $\| h \|_s < + \infty$, then it can be shown \cite{LermaLasserreMDP,AC_MDP_survey} that the stationary policy
\begin{align*}
	w_*(x) = \argmin_{u \in \sU(x)} \left\{ c(x,u) + P_u h(x) \right\}
\end{align*}
is optimal, and has average cost $\lambda$ for every $x$. The function $h$ is called the {\em relative value function}.

\subsection{Linearly solvable MDPs}
  
In a series of papers \cite{TodorovNIPS,TodorovCDC,TodorovPNAS},
Todorov has introduced a class of Markov decision processes, for which
solving the ACOE reduces to solving an eigenvalue problem. In this set-up, which we have described
informally in Section~\ref{ssec:setup_examples}, the action space
$\sU$ is the probability simplex $\cP(\sX)$, which is compact in the
Euclidean topology, and for each $u \in \cP(\sX)$ we have $P_u(x,y)
\deq u(y),\forall (x,y) \in \sX \times \sX$. Thus, any state feedback law (Markov policy)
$w: \sX \rightarrow \cP(\sX)$ induces the state transitions directly via
\begin{align*}
\Pr\lbrace X_{t+1} = y \vert X_t = x \rbrace = P_{w(x)}(x,y) = [w(x)](y), \quad t \ge 1.
\end{align*}
In other words, if $X_t = x$, then $u(\cdot) = w(x)$ is the probability
distribution of the next state $X_{t+1}$. Hence, there is a one-to-one
correspondence between Markov policies $w$ and Markov matrices
$P \in \cM(\sX)$, given by $w(x) = P(x,\cdot)$.

To specify an MDP, we fix a state cost function $f \in \cC_+(\sX)$ and a Markov matrix $\PD$ as the passive dynamics, which specifies the state transition probabilities in the absence of control. The one-step state-action cost function $c(x,u)$ is given by \eqref{eq:Todorov_cost}. If we  use  the shorthand  $c(x,  P)$  for $c(x,  P(x,\cdot))$, then the average cost  of a policy $P \in \cM(\sX)$ starting at $X_1 = x_1$ can be written as
\begin{align*}
J(P,x_1) = \limsup_{T \rightarrow \infty}  \frac{1}{T} \E_{x_1}^P \left[ \sum_{t=1}^T c(X_t, P) \right].
\end{align*} 
Intuitively, if $P$ has a small average cost, then the induced Markov chain $\bd{\mathit X} = \lbrace X_t \rbrace$ has a small average state cost, and its one-step transitions stay close to those prescribed by $\PD$.

The ACOE for this problem takes the form
\begin{align}\label{eq:ACOE_f}
h(x) + \lambda = f(x) + \min _{u \in \cP(\sX)} \lbrace D(u \| \PD(x,\cdot)) + \E_u h \rbrace.
\end{align}
For a given $h \in \cC(\sX)$, the minimization of the right-hand side of \eqref{eq:ACOE_f} can be done in closed form. To see this, let us define, for every $\varphi \in \cC(\sX)$, the {\em twisted kernel} \cite{BalajiMeyn}
\begin{align}\label{eq:twisted_kernel}
\twP_{\varphi}(x,\cdot) \deq \frac{\PD(x,\cdot)e^{-\varphi(\cdot)}} {\PD e^{-\varphi} (x)}, \qquad x \in \sX
\end{align}
which is obviously an element of $\cM(\sX)$. Then we have
\begin{align}
	 & \min_{u \in \cP(\sX)} \{  D(u \| \PD(x,\cdot)) + \E_u h \}  \nonumber \\
	 &=   \min_{u \in \cP(\sX)} \left\{ \E_u \left[ \log \frac {u(Y)}{\PD(x,Y)} + h(Y) \right] \right\} \nonumber \\
		&=  \min_{u \in \cP(\sX)} \left \{ \E_u \left[ \log \frac {u(Y)}{\twP_h(x,Y)} \right] - \log \PD e^{-h}(x) \right \} \label{eq:acoe_r_3}
\end{align}
If we further define $\Lambda_h(x) \deq \PD e^{-h}(x)$, then the quantity in braces in \eqref{eq:acoe_r_3} can be written as $D(u \| \twP_h(x,\cdot)) - \log \Lambda_h(x)$. Using the fact that the divergence $D(\mu \| \nu)$ between any two $\mu,\nu \in \cP(\sX)$ is nonnegative and equal to zero if and only if $\mu = \nu$ \cite{CoverThomas}, we see that the minimum value in \eqref{eq:acoe_r_3} is uniquely achieved by $u_*(x) = \twP_h(x,\cdot)$ and is equal to $-\log \Lambda_h (x)$.  Thus, we can rewrite the ACOE \eqref{eq:ACOE_f} as
\begin{align}\label{eq:acoe22}
h(x) + \lambda =  f(x) - \log \Lambda_h(x), \qquad \forall x \in \sX.
\end{align}
If we now consider the {\em exponentiated} relative value function $V \deq e^{-h}$, then \eqref{eq:acoe22} can be also written as $e^{-f}\PD V(x) = e^{-\lambda} V(x)$. Expressing this in vector form, we obtain the so-called {\em multiplicative Poisson equation} (MPE) \cite{BalajiMeyn}:
\begin{align}\label{eq:MPE_0}
e^{-f}\PD V = e^{-\lambda} V
\end{align}
To construct the optimal policy for our MDP, we first solve the MPE \eqref{eq:MPE_0} for $ \lambda$ and $V$, obtain $h$, and then compute the twisted kernel $\twP_{h}(x,\cdot)$ for every $x \in \sX$. The MPE is an instance of a so-called {\em Frobenius--Perron eigenvalue} (FPE) problem \cite{Seneta}; there exist efficient methods for solving such problems, e.g., a recent algorithm due to Chanchana \cite{Prakash}. We also should point out that, for each $x \in \sX$, the twisted kernel \eqref{eq:twisted_kernel} is a {\em Boltzmann--Gibbs distribution} on the state space $\sX$ with energy function $h$ and base measure $\PD(x,\cdot)$. Boltzmann--Gibbs distributions arise in various contexts, e.g., in statistical physics and in the theory of large deviations \cite{Ellis,Streater}, as solutions of variational problems over the space of probability measures that involve minimization of a Gibbs-type free energy functional, consisting of an affine ``energy" term and a convex ``entropy'' term (given by the divergence relative to the base measure). Indeed, the functional being minimized on the right-hand side of \eqref{eq:ACOE_f} is precisely of this form.

In the sequel, we will often need to consider simultaneously several MDPs with different state costs $f$. Thus, whenever need arises, we will indicate the dependence on $f$ using appropriate subscripts, as in $c_f, \lambda_f, h_f, V_f$, etc. For instance, the MPE \eqref{eq:MPE_0} for a given state cost $f$ is
\begin{align}\label{eq:MPE}
	\PD_f V_f = e^{-\lambda_f} V_f,
\end{align}
where $\PD_f \deq e^{-f}\PD$, i.e., $\PD_f(x,y) = e^{-f(x)}\PD(x,y)$ for all $x,y \in \sX$.

\subsection{Some properties of Todorov's optimal policy}
\label{ssec:Todorov_properties}

We now investigate the properties of Todorov's optimal policy under the assumptions on the passive dynamics $\PD$ that are listed in Section~\ref{ssec:main_result}. Most of the results of this section are new (with some exceptions, which we point out explicitly); the proofs are given in the Appendix.

We start with the following basic existence and uniqueness result, which is implicit in \cite{TodorovNIPS}:
	
\begin{proposition}\label{pps:MPE_sol} Under Assumption~\ref{as:irred}, for any state cost $f \in \cC_+(\sX)$ the MPE \eqref{eq:MPE} has a strictly positive solution $V_f \in \cC_+(\sX)$ with the associated strictly positive eigenvalue $e^{-\lambda_f}$, and the only nonnegative solutions of \eqref{eq:MPE} are positive multiples of $V_f$. Moreover, the corresponding twisted kernel $\twP_{h_f}$ is also irreducible and aperiodic, and has a unique invariant distribution $\twpi_f = \twpi_f \twP_f \in \cP(\sX)$.
\end{proposition}
\begin{IEEEproof} Appendix~\ref{app:proof_MPE_sol}.\end{IEEEproof}
	
\noindent Since $V_f = e^{-h_f}$, the fact that any positive multiple of $V_f$ is a solution of the MPE is equivalent to the well-known fact that the relative value function $h_f$ as a solution of the ACOE \eqref{eq:ACOE_f} is unique up to additive constants. That is, if a particular $h_f$ solves \eqref{eq:ACOE_f}, then so does any $h_f + c$ for any additive constant $c \in \Reals$. For this reason, we can fix an arbitrary $x^\circ \in \sX$ and assume that $h_f(x^\circ) = 0$ for any $f$. This ensures that the mapping
\begin{align}\label{eq:cost_to_h}
	f \longmapsto h_f, \qquad h_f(x^\circ) = 0 
\end{align}
is well-defined. The following results are new:

\begin{proposition}\label{pps:h_bounded} Under Assumption~\ref{as:irred}, the mapping \eqref{eq:cost_to_h} is bounded on compact subsets of the cone $\cC_+(\sX)$: for any $f \in \cC_+(\sX)$,
	\begin{align}\label{eq:hbound}
		\| h_f \|_s \le \log \theta^{-1} + \bar{n} \| f \|_\infty,
	\end{align}
where $\bar{n}$ and $\theta$ are defined in Section~\ref{ssec:main_result}. Hence,
\begin{align}\label{eq:uniform_hbound}
	\sup_{f \in \cC_+(\sX); \, \| f \|_\infty \le C} \| h_f \|_s \le \log \theta^{-1} + \bar{n} C.
\end{align}
\end{proposition}
\begin{IEEEproof}Appendix~\ref{app:proof_h_bounded}.\end{IEEEproof}

\noindent Moreover, the dependence of the relative value function $h_f$ on the state cost $f$ is {\em continuous}:
\begin{proposition}\label{pps:h_cont} Under Assumptions~\ref{as:irred} and \ref{as:Dobrushin}, the mapping \eqref{eq:cost_to_h} is Lipschitz-continuous on compact subsets of $\cC_+(\sX)$: for every $C > 0$ there exists a constant $K = K(C) > 0$, such that for any two $f,g \in \cC_+(\sX)$ with $\|f \|_\infty, \| g \|_\infty \le C$ we have
	\begin{align}\label{eq:h_cont}
		\| h_f - h_g \|_s \le K \| f - g \|_\infty.
	\end{align}
\end{proposition}
\begin{IEEEproof}Appendix~\ref{app:proof_h_cont}.\end{IEEEproof}
	
\noindent More generally, the twisted kernel $\twP_\varphi$ depends smoothly on the ``twisting function'' $\varphi$:
\begin{proposition}\label{pps:twisted_kernels} Fix any two functions $\varphi,\varphi' \in \cC(\sX)$. Then the twisted kernels \eqref{eq:twisted_kernel} have the following properties: for any $x \in \sX$,
		\begin{align}
			D(\twP_\varphi(x,\cdot) \| \twP_{\varphi'}(x,\cdot)) &\le \frac{1}{8} \| \varphi - \varphi' \|^2_s \label{eq:twisted_KL}\\
			\| \twP_{\varphi}(x,\cdot) - \twP_{\varphi'}(x,\cdot) \|_1 &\le \frac{1}{2} \| \varphi - \varphi' \|_s. \label{eq:twisted_TV}
		\end{align}
Moreover, if Assumptions~\ref{as:irred} and \ref{as:Dobrushin} hold, then there exists a mapping $\kappa : \Reals_+ \to [0,1)$, such that
\begin{align}\label{eq:twisted_alpha_bound}
\| \varphi \|_s \le C \quad \Longrightarrow \quad	\alpha(\twP_\varphi) \le \kappa(C).
\end{align}
\end{proposition}
\begin{IEEEproof}Appendix~\ref{app:proof_twisted_kernels}.\end{IEEEproof}
	
\noindent We close with the following basic but important result on steady-state optimality:
\begin{proposition}\label{pps:Jproof} For any $f \in \cC_+(\sX)$ and any $P \in \cM_1(\sX)$, define
	\begin{align*}
		\bar J_f(P) \deq \E_{\pi_P}[c_f(X,P)] \equiv \E_{\pi_P}[J_f(P,X)].
	\end{align*}
Then
\begin{align*}
\bar J_f(\twP_{h_f}) = \inf_{P \in \cM_1(\sX)} \bar J_f(P).
\end{align*}
\end{proposition}
\begin{IEEEproof} Appendix~\ref{app:Jproof}.\end{IEEEproof}

\section{The proposed strategy}
\label{sec:strategy}

Our construction of a Hannan-consistent strategy for the problem of Section~\ref{sec:problem} is similar to the approach of Yu et al. \cite{Yu}. The main idea behind it is as follows. We partition the set of time indices $1,2,\ldots$ into nonoverlapping contiguous segments (phases) of increasing duration and, during each phase, use Todorov's optimal policy matched to the average of the state cost functions revealed during the preceding phases. As in \cite{Yu}, the phases are sufficiently long to ensure convergence to the steady state within each phase, and yet are sufficiently short, so that the policies used during successive phases are reasonably close to one another.

The phases are indexed by $m \in \Naturals$, where we denote the $m$th phase by $\cT_m$ and its duration by $\tau_m$. Given $\epsilon \in (0,1/3)$, we let $\tau_m = \lceil m^{1/3-\epsilon} \rceil$. The parameter $\epsilon$ is needed to control the growth of the total length of each fixed number of phases relative to the length of the most recent phase (we comment upon this in more detail in the next section). We also define $\cT_{1:m} \deq \cT_1 \cup \ldots \cup \cT_m$ (the union of phases $1$ through $m$) and denote its duration by $\tau_{1:m}$. Given a sequence $\{f_t\}$ of state cost functions, we define for each $m$ the average state costs
\begin{align*}
	\wh{f}^{(m)} \deq \frac{1}{\tau_m}\sum_{t \in \cT_m}f_t, \qquad
	 \wh{f}^{(1:m)} \deq \frac{1}{\tau_{1:m}} \sum_{t \in \cT_{1:m}} f_t
\end{align*}
and let $\wh{f}^{(0)} = \wh{f}^{(1:0)} = 0$. Our strategy takes the following form:
\begin{center}
\begin{tabular}{|l|}
\hline
for $m = 1,2,\ldots$\\
\ \ solve the MPE $e^{-\wh{f}^{(1:m-1)}}\PD e^{-h^{(m)}} = e^{-\lambda^{(m)}}e^{-h^{(m)}}$\\
\ \ let $P^{(m)} = \twP_{h^{(m)}}$\\
\ \  for $t \in \cT_{m}$\\
\ \ \ \ draw $X_{t+1} \sim P^{(m)}(X_{t},\cdot)$\\
\ \ end for\\
end for\\
\hline 
\end{tabular}
\end{center}
Since we use the same policy throughout each phase, the evolution of the state induced by the above algorithm is described by the following inhomogeneous Markov chain:
\begin{align*}
X_1 \xrightarrow{P^{(1)}} X_2 \xrightarrow{P^{(1)}} \ldots \xrightarrow{P^{(1)}} X_{\tau_1} \xrightarrow{P^{(2)}} X_{\tau_1 +1} \xrightarrow{P^{(2)}} \ldots
\end{align*}
The implementation of this strategy reduces to solving a finite-dimensional Frobenius--Perron eigenvalue (FPE) problem \cite{Seneta} at the beginning of each phase to obtain a Todorov-type relative value function. The corresponding twisted kernel then determines the stationary policy to be followed throughout that phase. An efficient method for solving FPE problems was recently developed by Chanchana \cite{Prakash}. This method makes use of the well-known Collatz formula for the FPE \cite{Seneta} and Elsner's inverse iteration algorithm for computing the spectral radius of a nonnegative irreducible matrix \cite{Ludwig}. It is an iterative algorithm, which at each iteration performs an LU factorization of an $|\sX|\times|\sX|$ matrix. The time complexity of each iteration is $O(|\sX|^3)$. Chanchana's algorithm outperforms the three best known algorithms for solving FPE problems, which all rely on Elsner's inverse iteration and have quadratic convergence. Numerical experimental results can be found in \cite[Section~3.5]{Prakash}.

\section{Proof of Theorem~\ref{thm:main_r}}
\label{sec:proof}

\subsection{The main idea}

Following the general outline in \cite{Yu}, the proof of Theorem~\ref{thm:main_r} can be divided into four major steps. The first step is to show that there is no loss of generality in considering a different notion of regret, i.e., the {\em steady-state regret}, which is the difference between the cumulative cost of the proposed strategy and the steady-state cost of a fixed stationary policy.  The second step is to bound the difference between the expected total cost of our strategy and the sum of expected steady-state costs within each phase. That is, for each $m$, the steady-state expectation of the cost incurred in phase $m$ is taken w.r.t.\ the unique invariant distribution of $P^{(m)}$. After this step, we may only concentrate on expectations over invariant state distributions, which renders the problem much easier. For the third step, we show that the sum of steady-state expected costs is not much worse than what we would get if, at the start of each phase $m$, we also knew all the state cost functions to be revealed during phase $m$, i.e., if we used the ``clairvoyant'' policy $P^{(m+1)}$ in phase $m$. In the fourth step, we consider the sum of expected costs in each phase that could be attained if we knew all the state cost functions in advance and used the optimal policy w.r.t.\ the average of all the state cost functions throughout all the phases. We show that this expected cost is actually greater than the sum of expected costs of each phase when we only know the state cost functions one phase ahead. We then assemble the bounds obtained in these four steps to obtain the final bound on the regret of our strategy.

\subsection{Preliminary lemmas}

Before proceeding to the proof of Theorem~\ref{thm:main_r}, we present two lemmas that will be used throughout. The proofs of the lemmas rely heavily on the results of Section~\ref{ssec:Todorov_properties}, and are detailed in Appendices~\ref{app:proof_uniform_bounds3} and \ref{app:proof_Policy_cont}.

\begin{lemma}[Uniform bounds]\label{lem:uniform_bounds3} There exists constants $K_0 \ge 0$, $K_1 \ge 0$ and $0 \le \alpha < 1$, such that, for every $f \in \mathcal F$ and every $m \in \Naturals$,
\begin{align*}
\| c_f(\cdot,P^{(m)}) \|_\infty \le K_0, \quad \| h^{(m)} \|_s \le K_1, \quad \alpha (P^{(m)}) \le \alpha
\end{align*}
Moreover, the bound $\| c_f(\cdot,P) \|_\infty \le K_0$ holds for all $P \in \cM_1(\sX)$, such that $D(P(x,\cdot)\| \PD(x,\cdot)) < \infty$ for all $x \in \sX$.
\end{lemma}
\begin{lemma}[Policy continuity]\label{lem:Policy_cont} There exists a constant $K_2 \ge 0$, such that, for every $m \in \Naturals$,
\begin{align}\label{eq:kernel_continuity}
\| P^{(m+1)}(x, \cdot) - P^{(m)}(x,\cdot) \|_1 \le \frac{K_2  \tau_{m} }{\tau_{1:m}},
\end{align}
and
\begin{align}\label{eq:inv_dist_continuity}
\| \pi^{(m+1)} - \pi^{(m)} \|_1 \le \frac{K_2 \tau_{m}}{(1-\alpha) \tau_{1:m}}
\end{align}
where $\pi^{(m)}$ is the unique invariant distribution of $P^{(m)}$. Moreover, there exists a constant $K_3 \ge 0$, such that  for $D^{(m)}(x) \deq D(P^{(m)}(x,\cdot)\| \PD(x,\cdot)), \forall x \in \sX$, we have
\begin{align}\label{eq:Ddistance}
\| D^{(m)}(x) - D^{(m+1)}(x) \|_1 \le \frac{K_3 \tau_{m}}{\tau_{1:m}}.
\end{align}
\end{lemma}
\begin{remark} {\em As will be evident from the proof below, we can specify the  precise form of the regret bound in \eqref{eq:thm1} using the constants from the above lemmas:
\begin{align}\label{eq:thm2}
	& \sup_{P \in \cM_1^{\rho}(\sX)} \sup_{f_1, \ldots, f_T \in \cF} \frac{\E R_T(P)}{T} \nonumber\\
	&\le \frac{4}{3}\left(\frac{K_0(K_2 + 2)}{1-\alpha} + K_0 + K_3\right) T^{-1/4+\epsilon} + \frac{2K_0}{(1-\rho)T}.
\end{align}}
\end{remark}

\subsection{Details}

We are now ready to present the detailed proof of Theorem~\ref{thm:main_r}.

\noindent {\bf Step 1: Reduction to the steady-state case.} For any $P \in \cM_1(\sX)$, let us define the {\em steady-state regret} of our strategy $\gamma$ w.r.t.\ $P$ by
\begin{align*}
	R_T^{\text{ss}} (P) \deq C_T - \E_{\pi_P} \left[ \displaystyle \sum_{t=1}^{T} c_t(X, P) \right],
\end{align*}
which is the difference between the actual cumulative cost of $\gamma$ and the steady-state cost of the stationary unichain policy $P$ initialized with $\pi_P$. Now let us fix some $\rho \in [0,1)$ and consider an arbitrary $P \in \cM_1^{\rho}(\sX)$, where without loss of generality we can assume $D(P(x,\cdot)\| \PD(x,\cdot)) < \infty$ for all $x \in \sX$. For each $t \ge 1$, let $\nu_t = \delta_{x_1} P^{t-1}$ be the distribution of $X_t$ in the Markov chain induced by the transition matrix $P$ and initial state $X_1 = x_1$. For any $T$, we have
\begin{align}
	& \lvert R_T^{\text{ss}} (P) - R_T (P)\rvert \nonumber\\
&= \Bigg \lvert \E_{x_1}^P \left[ \sum_{t=1}^{T} c_t(X_t, P) \right] - \E_{\pi_P} \left[ \sum_{t=1}^{T} c_t(X, P) \right] \Bigg \rvert \nonumber \\
&= \Bigg \lvert \displaystyle \sum_{t=1}^{T}  \lbrace \E_{\nu_t} [c_t(X_t, P)] - \E_{\pi_P} [c_t(X, P)] \rbrace \Bigg \rvert  \nonumber \\
&\le \displaystyle \sum_{t=1}^{T} \| c_t(\cdot, P) \|_\infty \| \nu_t - \pi_P \|_1 \nonumber \\
&\le 2 K_0 \displaystyle \sum_{t=1}^{T} \rho^{t-1} \le \frac{2K_0}{1 - \rho}  \label{eq:step_1},
\end{align}
where the second inequality is by Lemma~\ref{lem:uniform_bounds3} and the fact that $P \in \cM_1^{\rho}(\sX)$. Therefore, it suffices to show that the bound in \eqref{eq:thm1} holds with $\E R_T^{\text{ss}} (P)$ in place of $\E R_T (P)$.

\noindent{\bf Step 2: Steady-state approximation within phases.} In this step, we approximate the cumulative cost within each phase by its steady-state value. Let $M$ denote the number of complete phases up to time $T$, i.e. $\tau_{1:M} \le T < \tau_{1:M+1}$ (simple algebra gives $M \le (4/3)T^{3/4+\epsilon}$). Then we can decompose the total cost as
\begin{align*}
	C_T &= \displaystyle \sum_{t=1}^{\tau_{1:M}} c_t(X_t, P_t) + \displaystyle \sum_{t=\tau_{1:M}+1}^T c_t(X_t, P_t)  \\
	&\le \displaystyle \sum_{t=1}^{\tau_{1:M}} c_t(X_t, P_t) + K_0 \tau_{M+1} = C_{\tau_{1:M}} + K_0 \tau_{M+1},
\end{align*}
where the inequality is by Lemma~\ref{lem:uniform_bounds3}. Since all state costs are nonnegative by hypothesis,
\begin{align*}
\E_{\pi_P} \left[ \sum_{t=1}^{T} c_t(X, P)  \right] \ge \E_{\pi_P} \left[ \sum_{t=1}^{\tau_{1:M}} c_t(X, P)  \right],
\end{align*}
which implies that
\begin{align}\label{eq:steadydiff}
R_T^{\rm ss} (P) \le R_{\tau_{1:M}}^{\rm ss} (P) + K_0 \tau_{M+1}.
\end{align}
For every time step $t$, let $\mu_t$ be the state distribution induced by our strategy when starting from initial state distribution $\mu_1 = \delta_{x_1}$. Note that the transition matrix at time $t$ is $P_t = P^{(m)}$ if $t \in {\mathcal T}_m$. We can decompose the expected cost in the first $M$ phases as
\begin{align}\label{eq:Mphasecost}
\E C_{\tau_{1:M}} = \displaystyle \sum_{m=1}^{M} \sum_{t \in {\mathcal T}_m} \E_{\mu_t} \left[c_t(X,P^{(m)})\right],
\end{align}
and for every $t \in {\mathcal T}_m$ we have
\begin{align*}
	&\E_{\mu_t} \left[c_t(X,P^{(m)})\right]  \\
	&\le \E_{\pi^{(m)}} \left[c_t(X,P^{(m)})\right] + \| c_t(\cdot, P^{(m)}) \|_\infty \| \mu_t - \pi^{(m)} \|_1 \\
	&\le \E_{\pi^{(m)}} \left[c_t(X,P^{(m)})\right] + K_0 \| \mu_t - \pi^{(m)} \|_1,
\end{align*}
where the last step is by Lemma~\ref{lem:uniform_bounds3}. In addition, for every $k \in \{0,1,\ldots, \tau_m-1\}$, we have
\begin{align*}
	&\| \mu_{\tau_{1:m-1}+k+1} - \pi^{(m)} \|_1 \\
	&= \left\| \mu_{\tau_{1:m-1}+1} {(P^{(m)})}^k- \pi^{(m)} {(P^{(m)})}^k  \right\|_1 \\
	&\le \alpha^k \| \mu_{\tau_{1:m-1}+1} - \pi^{(m)}  \|_1 \le 2 \alpha^k,
\end{align*}
where the first inequality is due to Lemma~\ref{lem:uniform_bounds3}. Hence,
\begin{align*}
	&\sum_{t \in {\mathcal T}_m} \E_{\mu_t} \left[c_t(X,P^{(m)})\right] \\
	&\le \sum_{t \in {\mathcal T}_m} \E_{\pi^{(m)}} \left[c_t(X,P^{(m)})\right] + 2K_0 \displaystyle \sum_{k=0}^{\tau_m -1} \alpha^k \\
	&\le \sum_{t \in {\mathcal T}_m} \E_{\pi^{(m)}} \left[c_t(X,P^{(m)})\right] + \frac {2K_0}{1-\alpha}.
\end{align*}
Substituting this into \eqref{eq:Mphasecost}, we have 
\begin{align*}
\E C_{\tau_{1:M}} \le \displaystyle \sum_{m=1}^{M} \sum_{t \in {\mathcal T}_m} \E_{\pi^{(m)}} \left[c_t(X,P^{(m)})\right] + \frac {2K_0 M}{1-\alpha}.
\end{align*}
\noindent{\bf Step 3: Looking one phase ahead.} In this step, we show that the steady-state cost in each phase is not much worse than what we could get if we knew everything one phase ahead. For every $m \in \{1,\ldots,M\}$, we have 
\begin{align*}
	&\sum_{t \in {\mathcal T}_m} \E_{\pi^{(m)}} \left[c_t(X,P^{(m)})\right] \\
	&\le \sum_{t \in {\mathcal T}_m} \E_{\pi^{(m+1)}} \left[c_t(X,P^{(m)})\right] + K_0 \tau_m \| \pi^{(m+1)} - \pi^{(m)} \|_1  \\
	&\le \tau_m \E_{\pi^{(m+1)}} \left[\wh f^{(m)} + D^{(m)}\right] + \frac{K_0 K_2 \tau_{m}^2}{(1-\alpha) \tau_{1:m}} \\
	&= \tau_m \bar J_{\wh f^{(m)}}(P^{(m+1)}) + \tau_m \E_{\pi^{(m+1)}} \left[D^{(m)} - D^{(m+1)}\right] \\
	&\qquad + \frac{K_0 K_2 \tau_{m}^2}{(1-\alpha) \tau_{1:m}} \\
	&\le \tau_m \bar J_{\wh f^{(m)}}(P^{(m+1)}) + \left(\frac{K_0 K_2}{1 - \alpha} + K_3\right) \frac{\tau_{m}^2}{\tau_{1:m}} ,
\end{align*}
where the first inequality is by Lemma~\ref{lem:uniform_bounds3}, the second inequality is by Lemma~\ref{lem:Policy_cont}, and the last inequality is due to \eqref{eq:Ddistance} in Lemma~\ref{lem:Policy_cont}. So we now have
\begin{align}
\E C_{\tau_{1:M}} &\le \displaystyle \sum_{m=1}^{M} \tau_m \bar J_{\wh f^{(m)}}(P^{(m+1)}) \nonumber \\
&\qquad + \displaystyle \sum_{m=1}^{M}\left(\frac{K_0 K_2}{1 - \alpha} + K_3\right) \frac{\tau_{m}^2}{\tau_{1:m}} + \frac {2K_0 M}{1-\alpha}. \label{eq:cumul_cost}
\end{align}
\noindent{\bf Step 4: Looking $M$ phases ahead.} In this step, we consider the fictitious situation where we know everything $M$ phases ahead, and show that the resulting steady-state value is actually greater than what we could get if we knew everything just one phase ahead. In other words, we claim that 
\begin{align}\label{eq:Jinduction}
\displaystyle \sum_{m=1}^{M} \tau_m \bar J_{\wh f^{(m)}}(P^{(m+1)}) \le \displaystyle \sum_{m=1}^{M} \tau_m \bar J_{\wh f^{(m)}}(P^{(M+1)}).
\end{align}
To see that this claim is true, we apply backward induction:
\begin{align*}
	&\displaystyle \sum_{m=1}^{M} \tau_m \bar J_{\wh f^{(m)}}(P^{(M+1)}) \\
	&= \displaystyle \sum_{m=1}^{M-1} \tau_m \bar J_{\wh f^{(m)}}(P^{(M+1)}) + \tau_M \bar J_{\wh f^{(M)}}(P^{(M+1)})  \\
	&= \tau_{1:M-1} \bar J_{\wh f^{(1:M-1)}}(P^{(M+1)}) + \tau_M \bar J_{\wh f^{(M)}}(P^{(M+1)})  \\
	&\ge \tau_{1:M-1} \bar J_{\wh f^{(1:M-1)}}(P^{(M)}) + \tau_M \bar J_{\wh f^{(M)}}(P^{(M+1)}) \\
	&= \displaystyle \sum_{m=1}^{M-1} \tau_m \bar J_{\wh f^{(m)}}(P^{(M)}) + \tau_M \bar J_{\wh f^{(M)}}(P^{(M+1)})  ,
\end{align*}
where the second equality is due to the fact that $\tau_{1:M-1} \wh f^{(1:M-1)} = \sum_{t \in {\mathcal T}_{1:M-1}} f_t  = \sum_{m=1}^{M} \tau_m \wh f^{(m)}$, while the inequality is by Proposition~\ref{pps:Jproof} and the fact that $P^{(M)} = \twP_{h^{(M)}}$, where $h^{(M)}$ is the relative value function for state cost $\wh f^{(1:M-1)}$. Repeating this argument, we obtain \eqref{eq:Jinduction}. Moreover,
\begin{align}
\sum_{m=1}^{M} \tau_m \bar J_{\wh f^{(m)}}(P^{(M+1)}) &=  \tau_{1:M} \bar J_{\wh f^{(1:M)}}(P^{(M+1)}) \nonumber\\
&  =  \tau_{1:M} \inf_{P \in \cM_1(\sX)} \bar J_{\wh f^{(1:M)}}(P) \nonumber\\
& = \inf_{P \in \cM_1(\sX)}  \E_{\pi_P} \left[ \displaystyle \sum_{t=1}^{\tau_{1:M}} c_t(X, P)  \right]. \label{eq:comparator}
\end{align}

After these four steps, we are finally in a position to bound the expected steady-state regret. Combining \eqref{eq:cumul_cost}--\eqref{eq:comparator}, we can write
\begin{align*}
	\E C_{\tau_{1:M}} &\le  \inf_{P \in \cM_1(\sX)}  \E_{\pi_P} \left[ \displaystyle \sum_{t=1}^{\tau_{1:M}} c_t(X, P)  \right] \\
	&\qquad + \sum_{m=1}^{M}\left(\frac{K_0 K_2}{1 - \alpha} + K_3\right) \frac{\tau_{m}^2}{\tau_{1:m}} + \frac {2K_0 M}{1-\alpha}.
\end{align*}
Therefore,
\begin{align}
	\E R_{\tau_{1:M}}^{\rm ss} (P) &=  \E C_{\tau_{1:M}} - \E_{\pi_P} \left[  \sum_{t=1}^{\tau_{1:M}} c_t(X, P)  \right]  \nonumber \\
	&\le \E C_{\tau_{1:M}} - \inf_{P \in \cM_1(\sX)} \E_{\pi_P} \left[  \sum_{t=1}^{\tau_{1:M}} c_t(X, P)  \right]  \nonumber \\
	&\le  \sum_{m=1}^{M}\left(\frac{K_0 K_2}{1 - \alpha} + K_3\right) \frac{\tau_{m}^2}{\tau_{1:m}} + \frac {2K_0 M}{1-\alpha}. \label{eq:phaseMregret}
\end{align}
Next we show that the right-hand side of \eqref{eq:phaseMregret} can be bounded by a quantity that is sublinear in $T$. From \eqref{eq:steadydiff}, we have
\begin{align*}
	\E R_T^{\rm ss} (P) &\le  \E R_{\tau_{1:M}}^{\rm ss} (P) + K_0 \tau_{M+1} \\
	&\le  \sum_{m=1}^{M}\left(\frac{K_0 K_2}{1 - \alpha} + K_3\right) \frac{\tau_{m}^2}{\tau_{1:m}} + \frac {2K_0 M}{1-\alpha} + K_0 \tau_{M+1}.
\end{align*}
Due to our construction of the phases, $M \le (4/3) T ^{3/4+\epsilon}$ and $\tau_{M+1} \le M$ if $M > 1$. Moreover, it is a matter of routine but tedious algebraic calculations to show that the choice $\tau_m = \lceil m^{\frac{1}{3} - \epsilon} \rceil$ for $m = 1, \ldots, M$ for any $\epsilon \in (0,1/3)$ is sufficient to guarantee that $\tau_m^2 \le \sqrt{\tau_{1:m}}$. Thus, we obtain
\begin{align*}
	\E R_T^{\rm ss} (P)& \le  \sum_{m=1}^{M}\left(\frac{K_0 K_2}{1 - \alpha} + K_3\right) \frac{\tau_{m}^2}{\tau_{1:m}} + \frac {2K_0 M}{1-\alpha} + K_0 M  \\
	&\le M\left(\frac{K_0(K_2 + 2)}{1-\alpha} + K_0 + K_3\right)  \\
&\le \frac{4}{3} \left(\frac{K_0(K_2 + 2)}{1-\alpha} + K_0 + K_3\right) T^{3/4+\epsilon}.
\end{align*}
Therefore, recalling \eqref{eq:step_1}, we finally obtain
\begin{align*}
	&\frac{\E R_T(P)}{T} \\
	&\le \frac{\E R_T^{\rm ss} (P)}{T} + \frac{2K_0}{T(1-\rho)} \\
	&\le \frac{4}{3} \left(\frac{K_0(K_2 + 2)}{1-\alpha} + K_0 + K_3\right)T^{-1/4+\epsilon} + \frac{2K_0}{T(1-\rho)},
\end{align*}
which completes the proof of Theorem~\ref{thm:main_r}.

\section{Simulations}
\label{sec:simulations}

In this section, we demonstrate the performance of our proposed
strategy on a simulated problem involving online (real-time) tracking
of a moving target on a large, connected, undirected graph $G$, which
models a terrain with obstacles. The state space is the set of all
vertices (nodes) of $G$. The target is executing a stationary random walk on $G$ with a randomly sampled transition probability matrix, which is different from the one that governs the passive dynamics $\PD$. The motion of both the tracking agent and the target must conform to the topology of $G$, in the sense that both can only move between neighboring vertices. The graph used in our simulation has $564$ vertices.

To make sure that Assumptions~\ref{as:irred} and \ref{as:Dobrushin} are satisfied, we construct the
passive dynamics in the form $\PD = (1-\delta) P_1 + \delta
P_0$ for some $\delta \in (0,1)$. Here, $P_1$ is a random walk that represents environmental
constraints, allowing the agent to go from a given node either to any
adjacent node (with equal probability) or to remain at the current
location. To ensure that the agent is sufficiently mobile, the
probability of not moving is chosen to be relatively small (in our case, $0.01$) compared
to the probability of transitioning to any of the neighboring
nodes. Since the underlying graph is connected, the random walk $P_1$
is irreducible; it is also aperiodic since $P_1(x,x) > 0$ for all vertices $x$.  We also add a perturbation random walk $P_0$, which has a
fixed column of ones (we can think of the node indexing that column as
a ``home base'' for the agent), and zeros elsewhere. The ``size'' of the perturbation is controlled by
$\delta$, which is set to be small (we have chosen $\delta = 0.01$), so the agent only has a
slight chance of returning to ``home base'' from any given node within
one step. This
perturbation ensures that no two rows of $\PD$ are orthogonal, and
$\alpha(\PD) \le 1-\delta = 0.99$. 
\begin{figure}[h]
\centering {
\includegraphics[scale = 0.40]{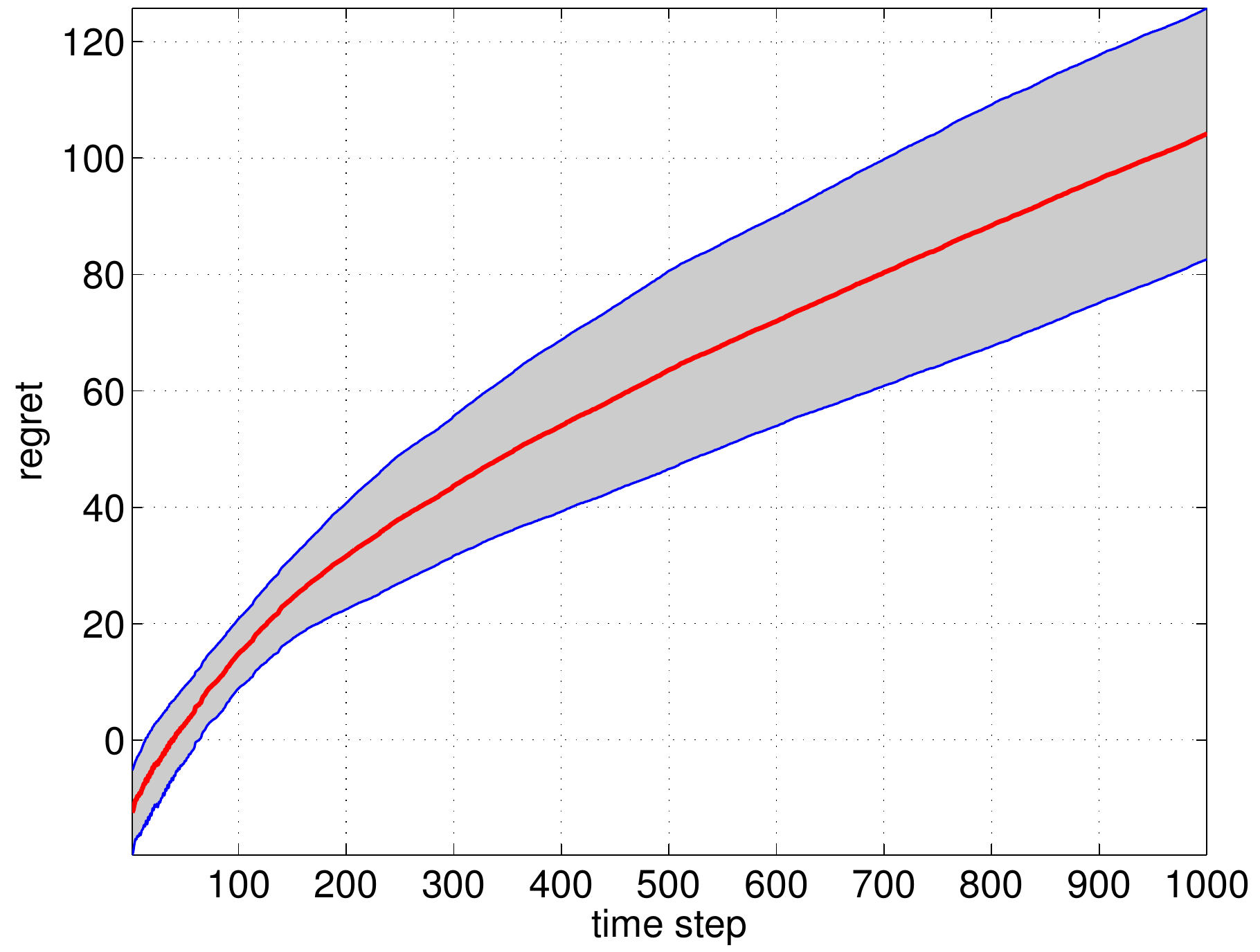}}
\caption{Regret versus time. The red curve shows the average of the regret (the difference between the total cost of our strategy up to each time $t$ and the total cost of the best stationary policy up to that time) over $100$ independent realizations of the simulation.  At each time $t$, the height of the gray region corresponds to one sample standard deviation.}
\label{fig:subfig1}
\end{figure}

The simulation consists of a number of independent experiments. Each individual experiment runs for $T = 1000$ time steps. We first randomly sample a transition matrix for the target motion. After simulating the target's random walk for $T$ steps, we record the target locations and use them to generate a sequence of state cost functions $\{ f_t\}^T_{t=1}$. Then we feed these $1000$ state cost functions sequentially to our online algorithm and compute the resulting cumulative cost $C_T$. At each time $t$, the tracking agent is in state (location) $x_t$, the target is at location $s_{t}$, and the agent's action is $P_t$. The cumulative cost after $T$ time steps is
\begin{align}
C_T =   \sum_{t=1}^T \left[f_t(x_t) + D\left(P_t(x_t,\cdot) \| \PD(x_t,\cdot)\right)\right], 
\end{align}
with state costs $f_t(x_t) = d_G (x_t, s_{t})$, where $d_G(\cdot,\cdot)$ is the graph distance (number of edges in the shortest path) between the agent's current location and the location of the target,
normalized by the diameter of $G$. Then we compute the best stationary policy $P$ in hindsight for the average of all the state costs by solving the MPE
\begin{align*}
e^{-\wh f}\PD e^{-h} = e^{-\lambda} e^{-h}
\end{align*}
for the relative value function $h$, where $\wh f = \frac{1}{T}\sum^T_{t =1}f_t$, and then setting
\begin{align*}
P(x,\cdot) = \frac{\PD(x,\cdot)e^{-h(\cdot)}} {\PD e^{-h} (x)}, \qquad x \in \sX
\end{align*}
 The regret is then computed with respect to the steady-state cost of this best stationary policy:
\begin{align*}
	R_T (P) = C_T - \E_{\pi_P} \left[ \displaystyle \sum_{t=1}^{T} c_t(X, P) \right],
\end{align*}
where
\begin{align*}
c_t(X,P) = f_t(X) + D\left(P(X,\cdot) \| \PD(X,\cdot)\right),
\end{align*} and $\pi_P$ is the unique invariant distribution of $P$.

To plot the regret versus time with error
  bars, we implement the experiment $100$ times and compute the
    empirical average of the regret across experiments. For
    each realization the agent was initialized with the same starting
  state. The evolution of the regret versus time is shown in
  Figure~\ref{fig:subfig1}, where the regret at time $t$ is defined as the
  total cost of our strategy up to time $t$ minus the total cost of
  the best stationary policy up to time $t$. We can see that the
  regret is growing sublinearly, as stated in Theorem~\ref{thm:main_r}. 

We also compare the total cost of our strategy to that of the best stationary baseline policy among a set $\tilde{{\mathcal N}}$ of $10^5$ randomly sampled stationary policies. Once again, each experiment runs for $T = 1000$ time steps. The baseline policy $P_{\rm baseline}$ is the one that has the smallest total cost
\begin{align*}
C_T(P_{\rm baseline}) =  \min_{P \in \tilde{{\mathcal N}}} \sum_{t=1}^T \left[f_t(x_t) + D\left(P(x_t,\cdot) \| \PD(x_t,\cdot)\right)\right].
\end{align*}
among the $10^5$ randomly sampled policies. The regret of our adaptive strategy is thus given by $C_T - C_T(P_{\rm baseline})$.

\begin{figure}[h]
\centering {
\includegraphics[scale = 0.40]{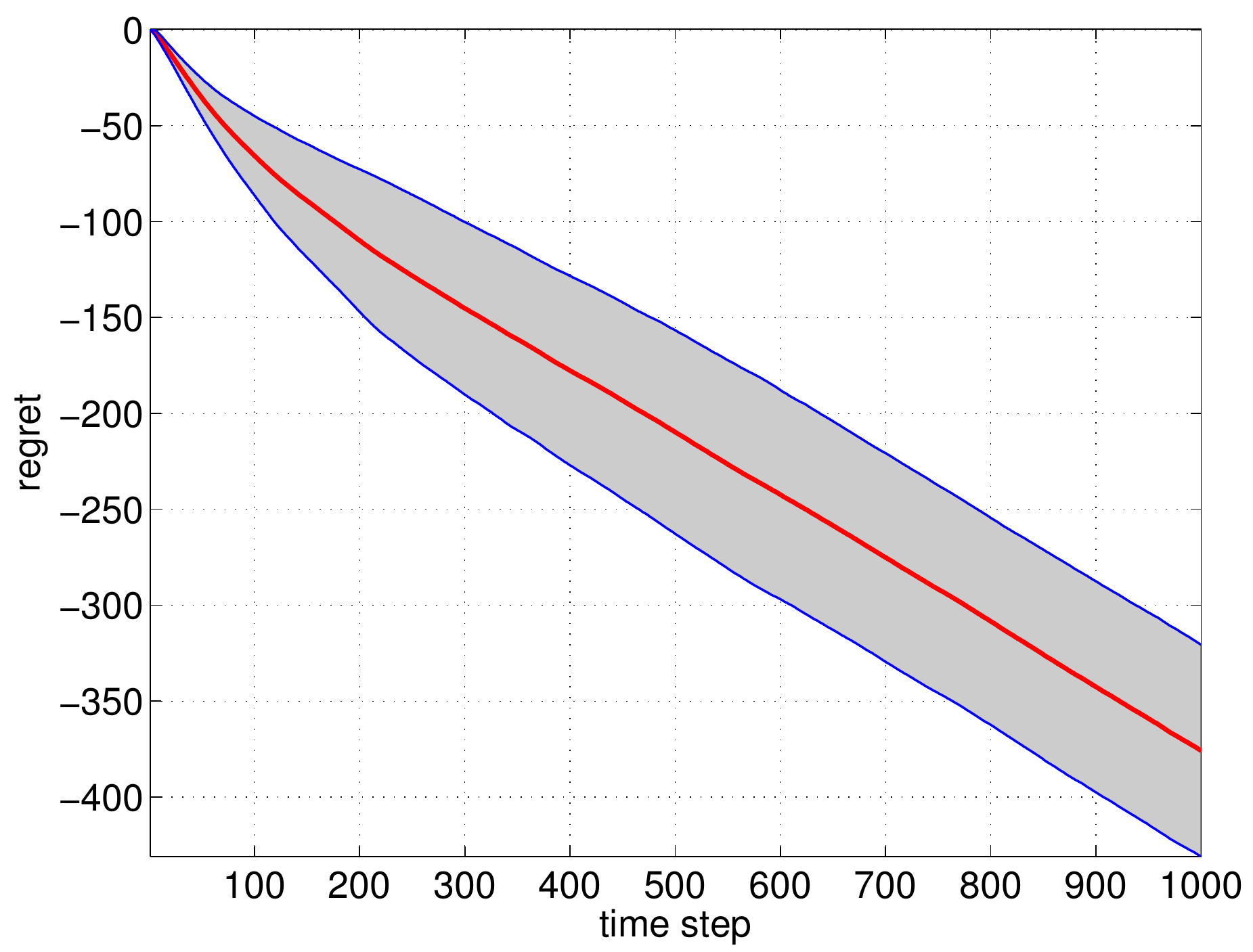}}
\caption{Comparison of our proposed strategy to the best stationary policy in a set of $10^5$ randomly sampled policies. The red curve shows the average of the regret (the difference between the total cost of our strategy up to each time $t$ and the total cost of the best stationary policy up to that time) over $100$ independent realizations of the simulation.  At each time $t$, the height of the gray area corresponds to one sample standard deviation.}
\label{fig:subfig}
\end{figure}

As before, there are $100$ independent experiments, where in each experiment  the agent using our strategy and the agent using the best sampled stationary policy were initialized with the same starting state. The evolution of the regret versus time is shown in Figure~\ref{fig:subfig}, where the regret at time $t$ is defined as the total cost of our strategy up to time $t$ minus the total cost of the best sampled stationary policy up to time $t$. We can see that the regret is negative, which implies that our strategy outperforms the best sampled stationary policy for each particular realization of the state cost sequence. 

\section{Conclusion and future work}

The problem studied in this paper combines aspects of both stochastic control and online learning. In particular, our construction of a Hannan-consistent strategy (a concept from the theory of online learning \cite{PLG}) uses several ideas and techniques from the theory of MDPs with average cost criterion, including some new results concerning optimal policies for  MDPs with KL control costs \cite{TodorovNIPS,TodorovCDC,TodorovPNAS}.

We have proved that, for any horizon $T$, our strategy achieves sublinear $O(T^{3/4})$ regret relative to any uniformly ergodic class of stationary policies, which is similar to the results of Yu et al. \cite{Yu} for online MDPs with finite state and action spaces. However, while our strategy (like that of \cite{Yu}) is computationally efficient, we believe that the $O(T^{3/4})$ scaling of regret with $T$ is suboptimal. Indeed, in the case when both the state and the action spaces are finite, Even-Dar et al. \cite{EvenDar} present a strategy that achieves a much better $O(\sqrt{T})$ regret. Of course, the strategy of \cite{EvenDar} involves recomputing the policy at {\em every} time step (rather than in phases, as is done here and in \cite{Yu}), which results in a significant loss of efficiency. An interesting open question, which we plan to address in our future work, is whether it is possible to attain $O(\sqrt{T})$ regret for online MDPs with KL control costs. A related challenge is to study these online MDPs in the (nonstochastic) bandit setting, where at each time step the agent only learns the value $f_t(X_t)$ of the state cost at time $t$ at the current state $X_t$, rather than the full state cost function $f_t \in \cC(\sX)$. While this bandit setting is more realistic, very little is known about it even for online MDPs with finite state and action spaces --- Neu et al.~\cite{onlineMDP_bandits}  constructed a strategy that achieves $O(T^{2/3}(\log T)^{1/3})$ regret, but it is not known whether this is optimal.

Another promising avenue for further research has to do with the apparent duality between our set-up and the theory of {\em risk-sensitive control} of Markov processes \cite{HHM,FHH1}. Indeed, the ACOE \eqref{eq:ACOE_f} can be viewed as a special case of the Isaacs equation for a certain dynamic two-player game with average cost criterion, in which Player 1 generates state cost functions, while Player 2 generates distributions over the state space (cf., e.g., \cite[p.~1805]{FHH1}). In the set-up of our Section~\ref{sec:problem}, Player 1 would correspond to the environment E, while Player 2 would be the agent A. We plan to explore this connection further.

Finally, as mentioned in the Introduction, we would like to extend our results to more general (e.g., compact) state spaces. This will require more sophisticated machinery, e.g., Foster--Lyapunov criteria and ergodicity w.r.t.\ weighted norms \cite{LermaLasserreMC,MeynTweedie}, as well as spectral theory of the MPE for Markov chains with general state spaces \cite{KontoMeyn}.

\appendix 

\subsection{Proof of Proposition~\ref{pps:MPE_sol}}
\label{app:proof_MPE_sol}

Consider the matrix $\PD_f \deq e^{-f}\PD$ with entries $\PD_f(x,y) = e^{-f(x)}\PD(x,y)$. For any $n \in \Naturals$ and any $x,y \in \sX$, $({\PD_f})^n(x,y) \ge e^{-n\| f \|_\infty} ({\PD})^n(x,y)$. Since $\PD$ is irreducible (Assumption~\ref{as:irred}), for any pair $x,y \in \sX$ of states there exists some $n \in \Naturals$, such that $({\PD})^n(x,y) > 0$. But then $({\PD_f})^n(x,y) > 0$ as well, which means that $\PD_f$ is also irreducible. Therefore, by the Frobenius--Perron theorem \cite{Seneta}, $\PD_f$ has a strictly positive right eigenvector $V_f$ with a positive eigenvalue $r$ (the Frobenius-Perron eigenvalue): $\PD_f V_f = r V_f$. Thus, $e^{-\lambda_f} = r$. Moreover, the FP eigenvalue is simple, and $\PD_f$ has no nonnegative right eigenvectors other than the positive multiples of $V_f$ \cite{Seneta}. This proves the existence and uniqueness part.

Now, using the fact that $V_f = e^{-h_f}$ solves the MPE \eqref{eq:MPE}, we can show that
\begin{align*}
	 \twP_{h_f}(x, y)
	&= e^{\lambda_f} \frac {V_f(y)}{V_f(x)}\PD_f(x,y),
\end{align*}
whence it follows that
\begin{align*}
(\twP_{h_f})^n(x, y) = e^{n\lambda^{*}_f} \frac {V_f(y)}{V_f(x)}({\PD_f})^n(x,y),
\end{align*}
As was just proved, $\PD_f$ is irreducible, and $V_f$ is strictly positive. Hence, for any pair $(x,y) \in \sX \times \sX$ there exists some $n \in \Naturals$, such that $(\twP_{h_f})^n(x, y) > 0$ as well. This proves the irreducibility of $\twP_{h_f}$. Now, since $\PD_f$ is irreducible, the Frobenius--Perron theorem says that there exists a unique strictly positive $\mu \in \cP(\sX)$, such that $\mu \PD_f = e^{-\lambda_f} \mu$ \cite{Seneta}. Now define $\twpi_f \in \cP(\sX)$ through
\begin{align*}
\twpi_f (x) \deq \frac{\mu (x) V_f(x)}{\sum_{y \in \sX}\mu(y)V_f(y)} \equiv \frac{\mu(x) V_f(x)}{\E_{\mu}V_f}, \qquad x \in \sX.
\end{align*}
A straightforward calculation shows that $\twpi_f$ is an invariant distribution of $\twP_{h_f}$. The uniqueness of $\twpi_f$ follows from the irreducibility of $\twP_{h_f}$.

\subsection{Proof of Proposition~\ref{pps:h_bounded}}
\label{app:proof_h_bounded}

We essentially follow the proof of Theorem~3.2 in \cite{FHH1}, with some simplifications. For each $T \in \Naturals$, define the function $W_T : \sX \to \Reals$ via
	$$
	e^{-W_T(x)} \deq \E_x \left[\exp\left(-\sum^{T}_{t=1}f(X_t)-h_f(X_{T+1})\right)\right],
	$$
	where $\E_x[\cdot]$ denotes the expectation w.r.t.\ the Markov chain $\bd{X} = (X_1,X_2,\ldots)$ with initial state $X_1 = x$ and transition matrix $\PD$. Then a simple inductive argument shows that
	\begin{align}\label{eq:W_ind}
	e^{-W_T(x)} = e^{-T\lambda_f-h_f(x)}.
\end{align}
	Indeed, for each $t$ let $\Psi_t \deq \prod^t_{s=1} \frac{V_f(X_s)}{\PD V_f(X_s)}$. Then, since $e^{-f(x)} = \frac{e^{-\lambda_f}V_f(x)}{\PD V_f(x)}$ by \eqref{eq:MPE}, we can write
	\begin{align}
		e^{-W_T(x)} &= e^{-T\lambda_f} \E_x \left[ \Psi_T V_f(X_{T+1})\right] \label{eq:W_ind_1} \\
		&= e^{-T\lambda_f} \E_x \left[\Psi_T \E[V_f(X_{T+1})|X_T]\right] \label{eq:W_ind_2}\\
		&= e^{-T\lambda_f} \E_x \left[\Psi_T \PD V_f(X_T)\right] \label{eq:W_ind_3}\\
		&= e^{-T\lambda_f} \E_x \left[ \Psi_{T-1} V_f(X_T)\right], \label{eq:W_ind_4}
	\end{align}
where \eqref{eq:W_ind_1} follows from definitions, \eqref{eq:W_ind_2} and \eqref{eq:W_ind_3} use the Markov property, and \eqref{eq:W_ind_4} again follows from definitions. Proceeding backwards, we get
\begin{align*}
	\E_x[\Psi_T V_f(X_{T+1})] &= \E_x[\Psi_1 V_f(X_2)] = \E_x [\Psi_1 \PD V_f(X_1)] \\
	&= V_f(x) = e^{-h_f(x)}.
\end{align*}
Substituting this into \eqref{eq:W_ind_1}, we get \eqref{eq:W_ind}, which in turn implies that $h_f(x) = W_T(x) - T\lambda_f$ for all $x \in \sX, T \in \Naturals$. Since $h_f(x^\circ) = 0$, we can write $h_f(x) = W_T(x) - W_T(x^\circ),  \forall x \in \sX, T \in \Naturals$. Let $\nu$ (respectively, $\nu^\circ$) be the distribution of $X_{\bar{n}+1}$ in the Markov chain with transition matrix $\PD$ and initial state $X_1 = x$ (respectively, $X_1 = x^\circ$). Then \begin{align}\label{eq:LR_lower}
\frac{\nu(y)}{\nu^\circ(y)} = \frac{(\PD)^{\bar{n}}(x,y)}{\nu^\circ(y)} \ge \theta > 0
\end{align}
for every $y \in \sX$. Consequently, for any $T > \bar{n}$ we have
\begin{align}
	&e^{-W_T(x)} \nonumber\\
	&= \E_x \left[ e^{-\sum^{\bar{n}}_{t=1} f(X_t)} e^{-\sum^T_{t=\bar{n}+1}f(X_t)-h_f(X_{T+1})}\right] \label{eq:W_diff_1}\\
	&\ge e^{-\bar{n} \| f \|_\infty} \E_x \left[ e^{-\sum^T_{t=\bar{n}+1} f(X_t)-h_f(X_{T+1})}\right] \label{eq:W_diff_2} \\
	&= e^{-\bar{n} \| f \|_\infty} \E_{x^\circ}\left[  e^{-\sum^T_{t=\bar{n}+1}f(X_t)-h_f(X_{T+1})} \frac{\nu(X_{\bar{n}+1})}{\nu^\circ(X_{\bar{n}+1})}\right] \label{eq:W_diff_3}\\
	&\ge \theta e^{-\bar{n} \| f \|_\infty} \E_{x^\circ} \left[  e^{-\sum^T_{t=1}f(X_t)-h_f(X_{T+1})} \right] \label{eq:W_diff_4}\\
	&= \theta e^{-\bar{n} \| f \|_\infty} e^{-W_T(x^\circ)}, \label{eq:W_diff_5}
\end{align}
where \eqref{eq:W_diff_1} is by definition, \eqref{eq:W_diff_3} follows from the Markov property and a change of measure, \eqref{eq:W_diff_4} follows from \eqref{eq:LR_lower} and from the fact that $f \ge 0$, and \eqref{eq:W_diff_5} is again by definition. Taking logarithms, we get $W_T(x) - W_T(x^\circ) \le \log \theta^{-1} + \bar{n}\|f \|_\infty, \forall T > \bar{n}$. Interchanging the roles of $x$ and $x^\circ$, we get $|h_f(x)|  \le \log \theta^{-1} + \bar{n}\|f \|_\infty$. This proves \eqref{eq:hbound}; \eqref{eq:uniform_hbound} follows immediately.

\subsection{Proof of Proposition~\ref{pps:h_cont}}
\label{app:proof_h_cont}

The basic idea is as follows. For a given $f \in \cC_+(\sX)$, let us introduce the dynamic programming operator $\DPop_f$ that maps any $\varphi \in \cC(\sX)$ to $\DPop_f \varphi \in \cC(\sX)$, where $\forall \varphi \in \cC(\sX), x \in \sX,$
	\begin{align*}
		\DPop_f \varphi (x) \deq f(x) + \inf_{\mu \in \cP(\sX)} \left\{ \E_\mu \varphi + D(\mu \| \PD(x,\cdot)) \right\}.
	\end{align*}
Then we can express the ACOE \eqref{eq:ACOE_f} as $h_f + \lambda_f = \DPop_f h_f$. Hence, for any $f,g \in \cC_+(\sX)$,
\begin{align}
	\left\| h_f - h_g \right\|_s &= \left\| \left( \DPop_f h_f - \lambda_f\right) -  \left( \DPop_g h_g - \lambda_g\right) \right\|_s \nonumber\\
	&= \left\| \DPop_f h_f - \DPop_g h_g \right\|_s \label{eq:span_bound_0}\\
	& \le \left\| \DPop_f h_f - \DPop_g h_f \right\|_s + \left\| \DPop_g h_f - \DPop_g h_g \right\|_s, \label{eq:span_bound}
\end{align}
where \eqref{eq:span_bound_0} uses the fact that the span seminorm is unchanged after adding a constant, and \eqref{eq:span_bound} is by the triangle inequality. We will then show the following:
\begin{enumerate}
	\item For any $\varphi \in \cC(\sX)$ and any $f,g \in \cC_+(\sX)$,
	\begin{align}\label{eq:DP_f_g}
		\left\| \DPop_f \varphi - \DPop_g \varphi \right\|_s \le 2 \| f - g \|_\infty.
	\end{align}
	\item For a fixed $f \in \cC_+(\sX)$, the dynamic programming operator $\DPop_f : \cC(\sX) \to \cC(\sX)$ is a contraction in the span seminorm: for every $M > 0$, there exists a constant $K' = K'(M) \in (0,1)$, such that for any two $\varphi,\varphi' \in \cC(\sX)$ with $\|\varphi \|_s, \| \varphi' \|_s \le M$ we have
	\begin{align}\label{eq:DP_contraction}
		\left\| \DPop_f \varphi - \DPop_f \varphi' \right\|_s \le K' \| \varphi - \varphi' \|_s.
	\end{align}
\end{enumerate}
Assuming that items 1) and 2) above are proved, we proceed as follows. First of all, the first term in \eqref{eq:span_bound} is bounded by $2 \| f - g \|_\infty$ by \eqref{eq:DP_f_g}. Next, since $\| f \|_\infty, \| g \|_\infty \le C$, Proposition~\ref{pps:h_bounded} guarantees that there exists some $M = M(C) < \infty$, such that $\| h_f \|_s, \| h_g \|_s \le M$. Therefore, there exists a constant $K' = K'(M) < 1$, such that the second term in \eqref{eq:span_bound} is bounded by $K' \| h_f - h_g \|_s$. Therefore, $\left\| h_f - h_g \right\|_s \le \frac{2}{1-K'} \| f - g \|_\infty$, which gives \eqref{eq:h_cont} with $K = 2/(1-K')$.

We now prove 1) and 2). For any function $\varphi \in \cC(\sX)$ and any two $f,g \in \cC_+(\sX)$, we have
\begin{align*}
	&\max_{x \in \sX} \left\{ \DPop_f \varphi(x) - \DPop_g \varphi(x)\right\} \\
&= \max_{x \in \sX}\Bigg\{\left[f(x) + \inf_{\mu \in \cP(\sX)}\left\{\E_\mu \varphi + D(\mu\|P^*(x,\cdot)) \right\}\right] \\
& \qquad \qquad - \left[g(x) + \inf_{\mu \in \cP(\sX)}\left\{\E_\mu \varphi + D(\mu\|P^*(x,\cdot)) \right\} \right]\Bigg\} \\
&= \max_{x \in \sX} \left[ f(x) - g(x)\right].
\end{align*}
Similarly, we get $\min_{x \in \sX} \left\{ \DPop_f \varphi(x) - \DPop_g \varphi(x)\right\}  = \min_{x \in \sX}\left[ f(x)-g(x)\right]$. Thus, $\| \DPop_f \varphi - \DPop_g \varphi \|_s = \| f - g \|_s \le 2 \| f - g \|_\infty$, so we have proved \eqref{eq:DP_f_g}.

To establish \eqref{eq:DP_contraction}, we follow the proof of Proposition~2.2 in \cite{MasiStettner} with some simplifications. Pick any $x,x' \in \sX$ and let
\begin{align*}
\nu = \argmin_{\mu \in \cP(\sX)} \left\{ \E_\mu\varphi' + D(\mu \| P^*(x,\cdot))\right\}, \\
\nu' = \argmin_{\mu \in \cP(\sX)} \left\{ \E_\mu\varphi + D(\mu \| P^*(x',\cdot))\right\},
\end{align*}
where explicitly $\nu(\cdot) = \twP_{\varphi'}(x,\cdot)$ and $\nu'(\cdot) = \twP_\varphi(x',\cdot)$.
Then
\begin{align*}
&	\left[ \DPop_f \varphi(x) - \DPop_f \varphi'(x)\right] - 	\left[ \DPop_f \varphi(x') - \DPop_f \varphi'(x')\right] \\
&= \inf_{\mu \in \cP(\sX)} \left\{\E_\mu \varphi + D(\mu \| P^*(x,\cdot)) \right\} \\
&\qquad - \inf_{\mu \in \cP(\sX)} \left\{ \E_\mu \varphi' + D(\mu \| P^*(x,\cdot))\right\} \\
&\qquad - \inf_{\mu \in \cP(\sX)} \left\{ \E_\mu \varphi + D(\mu \| P^*(x',\cdot)) \right\} \\
&\qquad + \inf_{\mu \in \cP(\sX)} \left\{ \E_\mu\varphi' + D(\mu \| P^*(x',\cdot))\right\} \\
&\le \E_\nu \varphi + D(\nu \| P^*(x,\cdot)) - \E_\nu \varphi' - D(\nu \| P^*(x,\cdot)) \\
&\qquad -\E_{\nu'}\varphi - D(\nu' \| P^*(x',\cdot)) + \E_{\nu'}\varphi' + D(\nu' \|  P^*(x',\cdot)) \\
&= \int (\varphi - \varphi') \mathrm{d}(\nu-\nu').
\end{align*}
A standard argument shows that that $\int (\varphi-\varphi') \mathrm{d}(\nu-\nu') \le \frac{1}{2} \| \varphi - \varphi' \|_s \| \nu - \nu' \|_1$. Consequently, 
\begin{align*}
	&\left\| \DPop_f \varphi - \DPop_f \varphi' \right\|_s \\
	&\le \frac{1}{2} \left\| \varphi - \varphi' \right\|_s \cdot \max_{x,x' \in \sX} \left\| \twP_{\varphi}(x,\cdot) - \twP_{\varphi'}(x',\cdot)\right\|_1.
\end{align*}
Then the proof of \eqref{eq:DP_contraction} will be complete if we can show that
\begin{align}\label{eq:DP_contraction_constant}
	K'(M) &\deq \frac{1}{2}\sup_{\varphi,\varphi'; \| \varphi \|_s, \| \varphi' \|_s \le M}\max_{x,x' \in \sX} \left\| \twP_{\varphi}(x,\cdot) - \twP_{\varphi'}(x',\cdot)\right\|_1 \nonumber \\
	&< 1.
\end{align}
Suppose that \eqref{eq:DP_contraction_constant} does not hold. Then there exist sequences $\{\varphi_n\}$, $\{\varphi'_n\}$ of functions with $\| \varphi_n \|_s, \| \varphi'_n \|_s \le M, \forall n$, a set $B \subset \sX$, and a pair of points $x,x' \in \sX$, such that
$$
\lim_{n \to \infty} \left[ \twP_{\varphi_n}(x,B) - \twP_{\varphi'_n}(x',B)\right] = 1,
$$
where for any $P \in \cM(\sX)$ we denote $P(x,B) \deq \sum_{y \in B}P(x,y)$. This implies in turn that
\begin{align}\label{eq:twP_limits}
\lim_{n \to \infty} \twP_{\varphi_n}(x,\sX \backslash B) =  \lim_{n \to \infty} \twP_{\varphi'_n}(x',B) = 0.
\end{align}
Since $\twP_\varphi(x,B) \ge e^{-\| \varphi \|_s} \PD(x,B)$, \eqref{eq:twP_limits} implies that $\PD(x,\sX \backslash B) =  \PD(x', B) = 0$. But this means that $\PD(x,B) - \PD(x',B) = 1$, which contradicts Assumption~\ref{as:Dobrushin}. Hence, \eqref{eq:DP_contraction_constant} holds.

\subsection{Proof of Proposition~\ref{pps:twisted_kernels}}
\label{app:proof_twisted_kernels}

We begin with \eqref{eq:twisted_KL}. From definitions, we have
	\begin{align}
		& D(\twP_{\varphi}(x,\cdot)\| \twP_{\varphi'}(x,\cdot)) \nonumber\\
		&\qquad = \E_{\twP_{\varphi}(x,\cdot)} [\varphi'(Y) - \varphi(Y)] + \log \frac{\Lambda_{\varphi'}(x)}{\Lambda_{\varphi}(x)}\label{eq:pps1_ineq}.
	\end{align}
A simple change-of-measure calculation shows that
	\begin{align}
		\frac{\Lambda_{\varphi'}(x)}{\Lambda_{\varphi}(x)} &= \frac{\sum_y \PD(x,y)e^{-\varphi'(y)}}{\Lambda_{\varphi}(x)} \nonumber \\
		&= \frac{\sum_y e^{\varphi(y)-\varphi'(y)}\PD(x,y)e^{-\varphi(y)}}{\Lambda_{\varphi}(x)} \nonumber \\
		&= \E_{\twP_{\varphi}(x,\cdot)} [e^{\varphi(Y)-\varphi'(Y)}]. \label{eq:lambda_ratio}
	\end{align}
To bound the right-hand side of \eqref{eq:lambda_ratio}, we recall the well-known {\em Hoeffding bound} \cite{Hoeffding}, which for our purposes can be stated as follows: For any $\mu \in \cP(\sX)$ and any $\psi \in \cC(\sX)$,
\begin{align*}
	\log \E_\mu e^{\psi} \le \E_\mu \psi + \frac{ \| \psi \|^2_s}{8}.
\end{align*}
Applying this bound gives
\begin{align*}
\log \frac{\Lambda_{\varphi'}(x)}{\Lambda_{\varphi}(x)} \le \E_{\twP_{\varphi}(x,\cdot)} [\varphi(Y) - \varphi'(Y)] + \frac{\| \varphi - \varphi' \|_s^2 }{8}.
\end{align*}
Substituting this bound into \eqref{eq:pps1_ineq}, we see that the terms involving the expectation of the difference $\varphi - \varphi'$ cancel, and we are left with \eqref{eq:twisted_KL}. To prove \eqref{eq:twisted_TV}, we use Pinsker's inequality, $\| P_1 - P_2 \|_1 \le \sqrt {2D(P_1 \| P_2)}$ \cite{CoverThomas}. To prove \eqref{eq:twisted_alpha_bound}, we follow essentially the same strategy as in the proof of Proposition~\ref{pps:h_cont} to show that $\kappa(C) \deq \sup_{\varphi; \| \varphi \|_s \le C} \alpha(\twP_\varphi)  < 1$ for every $C > 0$.

\subsection{Proof of Proposition~\ref{pps:Jproof}}
\label{app:Jproof}

Fix some $P \in \cM_1(\sX)$. If there exists some $x \in \sX$ such that $\pi_P(x) > 0$ and $D(P(x,\cdot) \| \PD(x, \cdot)) = + \infty$, then Proposition~\ref{pps:Jproof} holds trivially. Thus, there is no loss of generality if we assume that $D(P(x,\cdot) \| \PD(x, \cdot)) < + \infty, \forall x \in \sX$. Then
\begin{align*}
	 \bar J_f(P) &= \E_{\pi_P} [f(X) + D(P(X,\cdot) \| \PD(X, \cdot)) ] \nonumber\\
	&= \sum_{x \in \sX} \pi_P(x) \Bigg [ f(x) + \sum_{y \in \sX} P(x,y) \log \frac{P(x, y)}{\twP_{h_f}(x, y)} \nonumber \\
	&\qquad + \sum_{y \in \sX} P(x,y) \log \frac{\twP_{h_f}(x, y)}{\PD(x, y)}  \Bigg] \nonumber \\
	&= \sum_{x \in \sX} \pi_P(x) \Bigg [ f(x) + \E_{\pi_P} D(P(X,\cdot) \| \twP_{h_f}(X, \cdot)) \nonumber \\
	&\qquad + \sum_{y \in \sX} P(x,y) \log \frac{e^{-h_f(y)}}{\Lambda_{h_f}(x)}  \Bigg] \nonumber \\
	&\ge \sum_{x \in \sX} \pi_P(x) \Bigg [ f(x) + \sum_{y \in \sX} P(x,y) \log \frac{e^{-h_f(y)}}{\Lambda_{h_f}(x)}  \Bigg] \nonumber \\
	&= \E_{\pi_P} [f(X) -Ph_f(X) - \log \Lambda_{h_f}(X) ]  \nonumber \\
	&= \E_{\pi_P} [f(X) -h_f(X) - \log \Lambda_{h_f}(X) ] ,
\end{align*}
where the inequality is due to the fact that the KL divergence is always nonnegative, and the last step is due to the fact that $\pi_P$ is the invariant distribution of $P$. By the ACOE \eqref{eq:acoe22}, we know that $f(x) - h_f(x) - \log \Lambda_{h_f}(x) = \lambda_f$ for every $x \in \sX$. So we have $\bar J_f(P) \ge \lambda_f, \forall P \in \cM_1(\sX)$. Note that if we take the expectation $\E_{\twpi_f}[\cdot]$ of both sides of the ACOE \eqref{eq:ACOE_f}, we get
\begin{align*}
	&\E_{\twpi_f} [ h_f(X) + \lambda_f]  \\
	&= \E_{\twpi_f} [f(X) + D(\twP_{h_f}(X, \cdot) \| \PD(X, \cdot)) + \twP_{h_f} h_f(X) ]  \\
	&= \E_{\twpi_f}[f(X) + D(\twP_{h_f}(X, \cdot) \| \PD(X, \cdot))] + \E_{\twpi_f} [ h_f(X) ],
\end{align*}
where the last equality is due to the fact that $\twpi_f$ is the invariant distribution of $\twP_{h_f}$. Therefore, $\E_{\twpi_f}[f(X) + D(\twP_{h_f}(X, \cdot) \| \PD(X, \cdot))] =  \bar J_f(\twP_{h_f}) = \lambda_f$. So we now have $\bar J_f(P) \ge \bar J_f(\twP_{h_f})$ for any $P \in \cM_1(\sX)$, which completes the proof of Proposition~\ref{pps:Jproof}.

\subsection{Proof of Lemma~\ref{lem:uniform_bounds3}}
\label{app:proof_uniform_bounds3}

For every state $x \in \sX$, let $G_x$ denote the set of states that can be reached from $x$ in one step by the passive dynamics $\PD$, i.e., $G_x \deq \lbrace y : \PD(x,y) > 0 \rbrace$. Let us also define $p^*_x = \displaystyle \min_{y \in G_x} \PD(x,y)$ and $p^* = \displaystyle \min_{x \in \sX} p^*_x$. Since $P^{(m)} = \twP_{h^{(m)}}$, and $h^{(m)}$ is bounded by Proposition~\ref{pps:h_bounded}, we have $\supp(P^{(m)}(x, \cdot)) \subseteq \supp(\PD(x, \cdot)) \equiv G_x$. Therefore,
\begin{align*}
	&D(P^{(m)}(x, \cdot) \| \PD(x, \cdot)) = \sum_{y \in G_x} P^{(m)}(x, y)\log \frac{P^{(m)}(x, y)}{\PD(x, y)} \\
	&\le \log \frac{1}{p^*_x}, \qquad  \forall x \in \sX, m \in \Naturals
\end{align*}
and for any $f \in \mathcal F$
\begin{align*}
	&\| c_f(\cdot, P^{(m)})\|_\infty \\
	&\le \| f \|_\infty + \max_{x \in \sX}  D(P^{(m)}(x, \cdot) \| \PD(x, \cdot)) \le 1 + \log \frac{1}{p^*}.
\end{align*}
Thus, the first bound of Lemma~\ref{lem:uniform_bounds3} holds with $K_0 = 1 + \log \frac{1}{p^*}$. The same argument works for any $P \in \cM_1(\sX)$ that satisfies $D(P(x, \cdot) \| \PD(x, \cdot)) < \infty, \forall x \in \sX$. The second bound holds by Proposition~\ref{pps:h_bounded}, where $K_1 = \log \theta^{-1} + \bar{n}$. The third bound follows from the second bound, $\| h^{(m)} \|_s \le K_1$, and by Proposition~\ref{pps:twisted_kernels} with $\alpha = \kappa(K_1)$.

\subsection{Proof of Lemma~\ref{lem:Policy_cont}}
\label{app:proof_Policy_cont}

Let us recall that each $P^{(m)}$ is given by the twisted kernel $\twP_{h^{(m)}}$, where the relative value function $h^{(m)}$ arises from the solution of the MPE \eqref{eq:MPE} with state cost $\wh{f}^{(1:m-1)}$. Then
\begin{align}
	&\| P^{(m+1)}(x, \cdot) - P^{(m)}(x, \cdot) \|_1 \le \frac{1}{2} \| h^{(m+1)} - h^{(m)}\|_s \nonumber \\
	&\le \frac{K_2}{2} \| \hat f^{(1:m)} - \hat f^{(1:m-1)} \|_\infty \le \frac{K_2 \tau_{m}}{\tau_{1:m}} \label{eq:policy_cont2}
\end{align}
where the first step is by Proposition~\ref{pps:twisted_kernels}, the second by Proposition~\ref{pps:h_cont} with $K_2 = K(1)$, and the third by Lemma~4.3 in \cite{Yu}. This proves \eqref{eq:kernel_continuity}. Moreover, Proposition~\ref{pps:MPE_sol} guarantees that $P^{(m)} = \twP_{h^{(m)}}$ has a unique invariant distribution $\pi^{(m)}$. Therefore,
\begin{align*}
	&\| \pi^{(m)} - \pi^{(m+1)} \|_1 \\
	&= \| \pi^{(m)} P^{(m)} - \pi^{(m+1)} P^{(m+1)} \|_1\\
	&\le \| \pi^{(m)} P^{(m)} - \pi^{(m)} P^{(m+1)} \|_1 \\
	&\qquad +  \| \pi^{(m)} P^{(m+1)} - \pi^{(m+1)} P^{(m+1)} \|_1 \\
	&\le \| P^{(m)} - P^{(m+1)} \|_\infty + \alpha \| \pi^{(m)} - \pi^{(m+1)} \|_1 \\
	&\le \frac{K_2 \tau_{m}}{\tau_{1:m}} + \alpha \| \pi^{(m)} - \pi^{(m+1)} \|_1, 
\end{align*}
where the third inequality follows from \eqref{eq:policy_cont2}. Rearranging, we get \eqref{eq:inv_dist_continuity}.

Next, from the form of $P^{(m)}$ and $P^{(m+1)}$ we have
\begin{align}
	&D^{(m)}(x) - D^{(m+1)}(x) \nonumber \\
	&= \E_x^{(m+1)} [h^{(m+1)}] - \E_x^{(m)} [h^{(m)}] + \log \frac{\Lambda_{h^{(m+1)}} (x)}{\Lambda_{h^{(m)}} (x)}\label{eq:D_diffphase},
\end{align}
where $\E^{(m)}_x[\cdot]$ denotes expectation w.r.t.\ $P^{(m)}(x,\cdot)$, and we can follow the same steps we have used in \eqref{eq:lambda_ratio} to show that
\begin{align*}
	 \frac{\Lambda_{h^{(m+1)}} (x)}{\Lambda_{h^{(m)}} (x) }   =  \E_x^{(m)} \left[e^{h^{(m)}-h^{(m+1)}}\right].
\end{align*}
Using the Hoeffding bound \cite{Hoeffding}, we can write
\begin{align}\label{eq:lambda_hoff}
 \log \frac{\Lambda_{h^{(m+1)}} (x)}{\Lambda_{h^{(m)}} (x)} \le \E_x^{(m)} [h^{(m)}-h^{(m+1)}] + \frac{\| h^{(m)} - h^{(m+1)} \|_s^2 }{8}
\end{align}
Substituting \eqref{eq:lambda_hoff} into \eqref{eq:D_diffphase} and simplifying, we get
\begin{align*}
	&D^{(m)}(x) - D^{(m+1)}(x) \\
&\le  \E_x^{(m+1)} [h^{(m+1)}] - \E_x^{(m)} [h^{(m+1)}] + \frac{1}{8} \| h^{(m)} - h^{(m+1)} \|_s^2 \\
&\le \| h^{(m+1)} \|_s \cdot \| P^{(m)}(x, \cdot) - P^{(m+1)}(x, \cdot) \|_1 \nonumber \\
&\qquad + \frac{1}{8} \| h^{(m)} - h^{(m+1)} \|_s^2  \nonumber \\
&\le \frac{K_1 K_2 \tau_{m}}{\tau_{1:m}}+ \frac{1}{8} \| h^{(m)} - h^{(m+1)} \|_s^2  \nonumber \\
&\le \frac{K_1 K_2 \tau_{m}}{\tau_{1:m}}+ \frac{K_2^2 \tau_{m}^2 }{2 \tau_{1:m}^2}  \nonumber \\
&\le  \left(K_1 K_2 + \frac{K_2^2}{2}\right)\frac{\tau_{m}}{\tau_{1:m}}. \label{eq:D_diffphase3}
\end{align*}
Here, the third step uses the fact that $\| h^{(m)} \|_s \le K_1$ (Lemma~\ref{lem:Policy_cont}) and \eqref{eq:policy_cont2}, the fourth also uses \eqref{eq:policy_cont2}, and the last is due to the fact that $\frac{\tau_{m}}{\tau_{1:m}} < 1$. Letting $K_3 = K_1 K_2 + \frac{K_2^2}{2}$, we get \eqref{eq:Ddistance}.


\begin{biography}[{\includegraphics[width=1in,height=1.25in,clip,keepaspectratio]{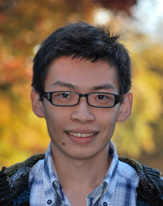}}]{Peng Guan} received the B.E.\ and M.Sc.\ degrees in Department of Automation from Tsinghua University, Beijing, China, in 2006 and 2009, respectively. He is currently working towards the Ph.D.\ degree in Electrical and Computer Engineering at Duke University. His research interests include stochastic control and online learning. 
\end{biography}

\begin{biography}[{\includegraphics[width=1in,height=1.25in,clip,keepaspectratio]{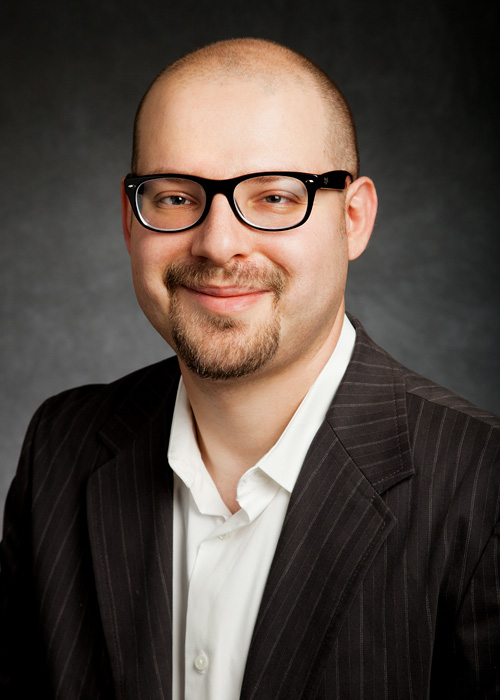}}]{Maxim Raginsky} received the B.S.\ and M.S.\ degrees in 2000 and the Ph.D.\ degree in 2002 from Northwestern University, Evanston, IL, all in electrical engineering. He has held research positions with Northwestern, the University of Illinois at Urbana-Champaign (where he was a Beckman Foundation Fellow from 2004 to 2007), and Duke University. In 2012, he has returned to UIUC, where he is currently an Assistant Professor with the Department of Electrical and Computer Engineering and the Coordinated
Science Laboratory. In 2013, Prof. Raginsky has received a Faculty Early Career Development (CAREER) Award from the National Science Foundation. His research interests lie at the intersection of information theory, machine learning, and control.
\end{biography}

\begin{biography}[{\includegraphics[width=1in,height=1.25in,clip,keepaspectratio]{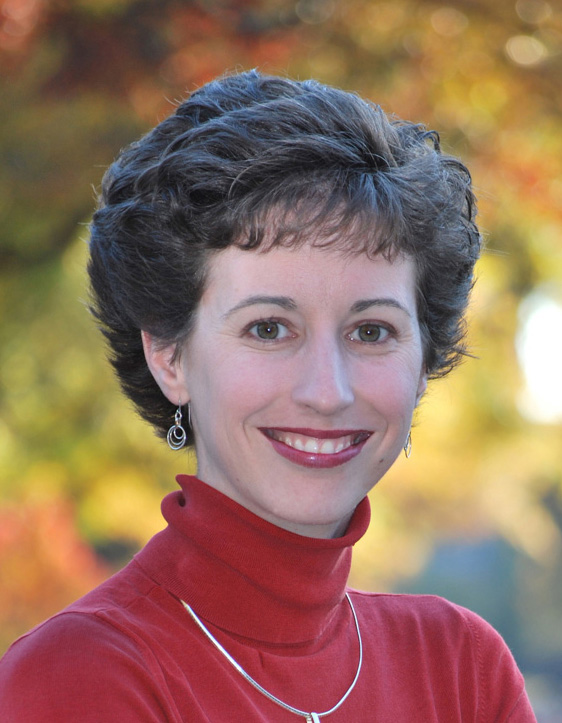}}]{Rebecca Willett} is an Associate Professor of Electrical and Computer Engineering at the University of Wisconsin-Madison. She completed her
PhD in Electrical and Computer Engineering at Rice University in 2005
and was an Assistant then Associate Professor of Electrical and
Computer Engineering at Duke University from 2005 to 2013. Willett received the National Science Foundation CAREER Award in 2007, is a
member of the DARPA Computer Science Study Group, and received an Air
Force Office of Scientific Research Young Investigator Program award
in 2010. Willett has also held visiting researcher positions at the
Institute for Pure and Applied Mathematics at UCLA in 2004, the
University of Wisconsin-Madison 2003-2005, the French National
Institute for Research in Computer Science and Control (INRIA) in
2003, and the Applied Science Research and Development Laboratory at
GE Healthcare in 2002. Her research interests include network and
imaging science with applications in medical imaging, neural coding,
astronomy, and social networks. Additional information, including
publications and software, are available online at
http://willett.ece.wisc.edu.
\end{biography}

\end{document}